%% file: main.tex
\documentclass[preprint,11pt,3p]{elsarticle}

\usepackage{booktabs} 
\usepackage{amssymb}
\usepackage{floatrow}
\usepackage{pgf-pie}
\usepackage{verbatim}
\usepackage[mathcal]{euscript}
\usepackage[final]{pdfpages}
\usepackage{pdfpages}
\usepackage[displaymath, mathlines]{lineno}
\usepackage{varioref} 
\usepackage{float}
\usepackage{amsfonts,amsmath}
\usepackage{tikz}
\usepackage{epsfig}
\usepackage{epstopdf}
\usetikzlibrary{shadings}
\usetikzlibrary{arrows}
\usetikzlibrary{shapes,snakes}
\usepackage{bm}
\usepackage{mathtools}
\usepackage{graphicx}
\usepackage{caption}
\usepackage{hyperref}
\hypersetup{colorlinks=true}
\usepackage[final]{pdfpages}
\usepackage{multicol}
\usepackage{tikz-3dplot} 
\usepackage{pdfpages}       
\usepackage{varioref} 
\usepackage{float}
\usepackage{multirow}

\usetikzlibrary{quotes,arrows.meta}
\numberwithin{equation}{section}  
\usepackage{tabularx}
\usepackage{arydshln} 
\usepackage{subfigure}
\usepackage{appendix} 
\usepackage{empheq}

\newtheorem{definition}{Definition}[section]

\newtheorem{remark}[definition]{Remark}

\newcommand \bei {\begin{itemize}}
\newcommand \eei {\end{itemize}}

\newcommand \del \partial

\begin{document}
\begin{frontmatter}
 \title{A semi-implicit second-order temporal scheme for solving the pressure head-based form of Richards’ and advection-dispersion equations}
 
\author{Nour-eddine Toutlini$^{a,b}$, Abdelaziz Beljadid$^{a,c,*}$, Azzeddine Soula\"imani$^b$}
  \address{$^a$ University Mohammed VI Polytechnic, Morocco}
 \address{$^b$ École de technologie supérieure, Canada}  
 \address{$^c$ University of Ottawa, Canada}

\cortext[mycorrespondingauthor]{Corresponding author.  
\newline
E-mail adresses: Abdelaziz.BELJADID@um6p.ma, abeljadi@uottawa.ca}

\begin{abstract}
In this study, a novel semi-implicit second-order temporal scheme combined with the finite element method for space discretization is proposed to solve the coupled system of infiltration and solute transport in unsaturated porous media. The Richards equation is used to describe unsaturated flow, while the advection-dispersion equation (ADE) is used for modeling solute transport.  The proposed approach is used to linearize the system of equations in time, eliminating the need of iterative processes. A free parameter is introduced to ensure the stability of the scheme. Numerical tests are conducted to analyze the accuracy of the proposed method in comparison with a family of second-order iterative schemes. The proposed numerical technique based on the optimal free parameter is accurate and performs better in terms of efficiency since it offers a considerable gain in computational time compared to the other methods. For reliability and effectiveness evaluation of the developed semi-implicit scheme, four showcase scenarios are used. The first two numerical tests focus on modeling water flow in heterogeneous soil and transient flow in variably saturated zones. The last numerical tests are carried out to simulate the salt and nitrate transport through unsaturated soils. The simulation results are compared with reference solutions and laboratory data, and demonstrate the effectiveness of the proposed scheme in simulating infiltration and solute transport through unsaturated soils.
\end{abstract}


\end{frontmatter}

\setcounter{tocdepth}{3}  
%

\section{Introduction}\label{S1}
\input{Docts/Intro}
\section{Governing equations}\label{S2}
\subsection{Richards equation}
\input{Docts/ARichards_equation}

\input{Docts/Ph.RHC}

\subsection{Advection-dispersion equation }
\input{Docts/ADE}

\section{Material and methods}\label{S3}
\subsection{Semi-discrete formulation}
\input{Docts/weak_form}
\subsection{Time discretization and proposed methods}
\input{Docts/timediscrete}

\subsection{Linearization technique}
\input{Docts/linearisation}
\section{Numerical results}\label{S4}
\begin{figure}
\centering
\includegraphics[scale=0.8]{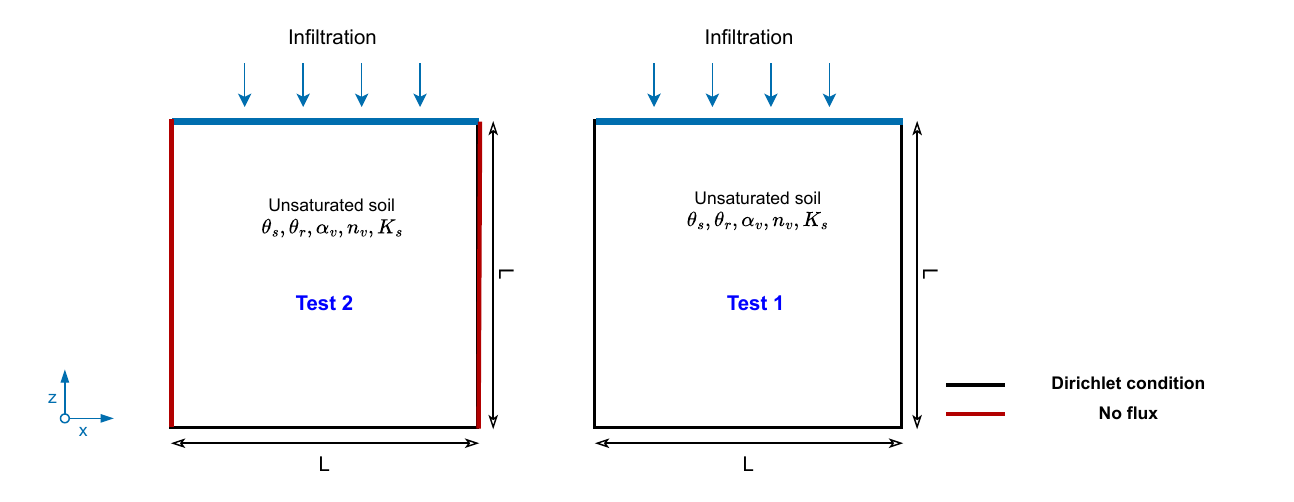} 
\caption{Schematic illustration of the two-dimensional infiltration  tests.}
\label{Fig:2dGA}
\end{figure}
In this section, we conduct numerical experiments in two-dimensional scenarios, considering both homogeneous and heterogeneous media and employing both the Gardner and the VG models for capillary pressure. 
Firstly, we will start with a numerical analysis of the studied methods by conducting numerically the convergence rate of the methods  \eqref{BDF2}-\eqref{SILf}. For the convergence analysis, we will use the $L^2$ error as described in \eqref{L2Error} to calculate the error that occurs between the reference or exact solution and the numerical solution. The rate of convergence, $p$,  will be computed using the following formula:
\begin{equation}
    p=\ln\left(\dfrac{||u -u_h||_{L^2} }{||u -u_{h/k}||_{L^2}}\right)\Big/ \ln\left(k\right),
    \label{order_of_convergence}
\end{equation}
where $u$ represents the reference solution or exact solution, $u_h$ represents the numerical solution and \(k\) represents the refinement factor between consecutive element sizes or time steps. The value \(k=2\) was selected in our numerical tests. Another criterion to test the efficiency of the methods is the CPU time (total runtime). We compare the numerical efficiency of the methods by looking at the CPU time required to reach a desired level of accuracy. 
Secondly, some application tests will be conducted for the efficient method. We should emphasize that the SILF2 method will be used with the optimal stabilized parameter \(\nu = 1\) as will be confirmed in the following section. Throughout our analysis,  all the numerical tests are performed using an 	Intel(R) Core(TM) i5-10310U CPU @ 1.70GHz, 2208 Mhz, 4 Core(s), 8 Logical Processor(s) using the finite element platform FreeFem++ \cite{Hecht2012, MR3043640} with the direct solver UMFPACK \cite{Davis2004}.
\subsection{Numerical tests - Convergence}
Here, we focus on the Richards equation as the first step since it is a nonlinear equation and the analysis will be effective to compare the presented methods. To analyze the convergence and efficiency of the proposed methods, we apply them to a variety of infiltration problems and compare them to reference solutions. We consider the following test problems:
\begin{enumerate}
\item[$\bullet$] Two infiltration tests with the Gardner model in homogeneous porous media for which analytical solutions exist \cite{tarcy2011}.  \item[$\bullet$] One infiltration test with the VG model in heterogeneous porous media for which the reference solution is used.
\end{enumerate}

\subsubsection{ Green and Ampt infiltration problem}
\input{Docts/numerical_tests}

\subsubsection{Heterogeneous soil test}
\input{Docts/Heter}
\newpage

In the following, we will focus on application tests related to water flow in variably saturated soil and solute transport through unsaturated soils. Specifically, we will consider:
\begin{enumerate}
\item[$\bullet$] A 2D numerical test modeling water flow in variably saturated sandy soil \cite{vauclin1979experimental}.
\item[$\bullet$] A 2D numerical test modeling water flow and soil salt transport in loamy soil \cite{haruzi2023modeling}.
\item[$\bullet$] A 2D numerical test modeling nitrate transport under surface drip fertigation \cite{li2005modeling}.
\end{enumerate}
\subsection{: 2D water table recharge experiment}
\input{Docts/Vaclintest}
\subsection{Soil salt transport}
\input{Docts/SolutePENNS}

\subsection{Nitrate transport}
\input{Docts/Nitrate}
\section{Conclusion}\label{S5}
A semi-implicit second-order temporal scheme (SILF2) is developed for modeling infiltration and solute transport in unsaturated soils. The pressure head-based form of the Richards equation is employed to express the infiltration process, whereas solute transport is described using the advection-dispersion equation. The finite element method is employed for spatial discretization. The developed SILF2 scheme includes a free stabilized parameter for which the optimal value is found to be 1.
For comparison, a family of second-order temporal schemes are employed, which require an iterative process to solve the Richards equation. The SILF2 scheme was evaluated using several test cases with different soils and boundary conditions. The results showed that the proposed approach based on the optimal free parameter \(\nu = 1\) is accurate and yields better results than the other methods in terms of computational cost. The proposed SILF2 scheme demonstrates its capabilities in simulating infiltration in scenarios involving heterogeneous soils and variably saturated zones. The reliability of the scheme is validated against reference solution and experimental data. The SILF2 scheme efficiently simulates solute transport in unsaturated soils, where two numerical tests are performed. One involving infiltration and soil salt transport in a homogeneous loamy soil, and the other focusing on nitrate transport under surface drip fertigation, which was compared with laboratory experimental data. The results highlight the SILF2 scheme's robust capabilities, making it a valuable tool for modeling infiltration and solute transport processes in soils.
\subsection*{\bf Acknowledgments}
Funding for this research was provided by UM6P/OCP Group of Morocco, the Moroccan Ministry of Higher Education, Scientific Research and Innovation and the OCP Foundation (APRD research program), and Natural Sciences and Engineering Research Council of Canada.

All the numerical tests presented in this study were conducted using the supercomputer simlab-cluster, supported by University Mohammed VI Polytechnic, and facilities of simlab-cluster HPC $\&$ IA platform.
\appendix
\renewcommand\thesection{Appendix~\Alph{section}}
\section{Truncation error}\label{appA:details}
\input{Docts/TError}

\newpage
\bibliographystyle{apalike}
\bibliography{ref}

\end{document}

%% file: Docts/Intro.tex
Modeling of infiltration and solute transport through unsaturated soils is of primary importance in tremendous
applications, including those relevant to the optimization of irrigation systems, hydrogeology, groundwater management, and environmental protection \cite{ bear2013dynamics, DAMODHARARAO2006, siyal2012minimizing, siyal2013strategies}. The prediction of these processes is complex and it can include processes such as evapotranspiration, plant water uptake, soil salinity, chemical reactions, and soil heterogeneity  \cite{Simunek2009}. Physical-based hydrological models have been developed to investigate the interactions of water and solutes in soils, which are conceptualized in the form of partial differential equations (PDEs). Various models for predicting the infiltration process are used \cite{beljadid2020continuum,  philip1969,Richards1931}. 
In our study, the traditional Richards equation \cite{Richards1931} is used to describe infiltration through unsaturated soils \cite{ Boujoudar2023,  Celia1990, clement2021, Keita2021, radu2016}. Conversely, the advection-dispersion equation \cite{freeze1979groundwater} is used to model the solute transport in soil \cite{li2005modeling, Simunek2016}. 
Analytical solutions for the Richards equation are restricted to simplified cases.  This pertains to the advection-dispersion equation since the solutions for the water content and water fluxes are required. Solving these PDEs numerically presents challenges in terms of stability, accuracy, and mass conservation, particularly in various heterogeneous, dry, and wet soil conditions. As a result, designing state-of-the-art numerical methods is required.
Several numerical methods have been developed to solve the Richards equation, such as finite element \cite{Celia1990, Keita2021}, finite volume \cite{eymard1999, Manzini2004}, and meshless \cite{Boujoudar2021, Boujoudar2023, Boujoudar2024Implicit}  methods. Similarly, the solute transport equation is addressed through various numerical techniques such as finite element \cite{ELAMRANI2023184, kaluarachchi1988finite, younes2022robust} and 
finite difference \cite{selim1981modeling} methods. The coupled system of infiltration and solute transport in unsaturated soils has been extensively examined concerning spatial discretization in numerous studies \cite{DIAW2001197, russo2012numerical, soraganvi2009modeling}. For example, in \cite{DIAW2001197}, the authors developed an efficient one-dimensional approach utilizing a combination of numerical methods: the finite element method for Richard's equation, a discontinuous finite element method for the advective term, and a finite difference method for the dispersive term. Some software programs are developed to solve this coupled system in porous media including Hydrus 1D, (2D/3D) \cite{ Simunek2005, simunek2006, Simunek2016}, HydroGeoSphere \cite{brunner2012}, etc. For more details about these software packages and the numerical methods that they used, one can refer to \cite{farthing2017, zha2019review}. Recently, a fresh perspective has come to light in tackling the water flow and solute transport in unsaturated soils, which is characterized by the integration of the physics-informed neural networks approach \cite{haruzi2023modeling}.  In this present study, we will focus on the finite element approach, which is widely used, because of its ability to deal with stiff problems \cite{ EYMARD2015186, girault2012,pinder2013}.
In the context of the infiltration, many classes of finite element methods have been applied to solve the Richards equation, such as the continuous Galerkin finite element \cite{Celia1990}, the mixed finite element \cite{Keita2021}, and the discontinuous Galerkin finite element methods \cite{antonietti2019discontinuous, clement2021}.

Regarding the time discretization techniques, implicit schemes are mostly used for solving the coupled system due to their broad stability range and enhanced accuracy \cite{radu2016, simunek2006}. Among the implicit schemes commonly employed for stiff problems, a notable choice is the family of backward differentiation formulas (BDFs). This category comprises a variety of methods based on one or multiple steps. These methods span from first order to fifth order, denoted as BDF1, BDF2, and so forth. It's worth noting that BDFs of higher orders (such as 6 and 7) tend to be avoided due to their inherent unconditional instability \cite{Durran}. In this study, we will focus on implicit second-order time-stepping techniques, including BDF2 and other multistep methods. The BDF2 time-stepping procedure has been popular in engineering literature due to its L-stable property and second-order accuracy \cite{Durran, Hundsdorfer,lambert1991numerical}. Some previous studies used the implicit multistep methods to solve the Richards equation, as stated in \cite{clement2021, Keita2021} and references therein. The main challenge of the multistep methods is they are not self-starting, which signifies that the solution at the first time step needs to be obtained by some other procedure. Conventionally, the initial time step solution is often achieved using the BDF1 method. However, some authors provide alternatives to these initialization methods to ensure more accuracy and stability \cite{nishikawa2019}. 
One notable challenge of multistep methods, in addition to lacking self-start capability, is dealing with the high nonlinear terms presented in the Richards equation, which necessitate an iterative method at each time step as we will discuss in the linearization part of this study.

Semi-implicit methods are an attractive alternative for solving the Richards equation. For instance, Keita et al. \cite{Keita2021} developed a linear semi-implicit second-order time stepping method with a mixed finite element method and they showed the advantages of this semi-implicit technique for solving the mixed form of the Richards equation. This motivates us to develop a noniterative semi-implicit method for solving the coupled system of unsaturated flow and solute transport in porous media. In the literature, to the best of our knowledge, there are few research studies on the implicit and semi-implicit methods addressing this coupled system \cite{kaluarachchi1988finite, ross2003modeling, srivastava1992three}. The main contribution of this study is the development of a noniterative, semi-implicit, second-order temporal scheme for solving the Richards equation in its pressure head-based form, alongside the advection-dispersion equation. A free parameter is introduced to obtain a stable semi-implicit scheme. This study includes a comparative analysis of the proposed approach against a family of iterative second-order temporal schemes, focusing on evaluating accuracy, stability, and computational cost. Numerical tests demonstrate that the semi-implicit scheme with the optimal free parameter yields more accurate solutions compared to other schemes. Additionally, we aim to demonstrate the effectiveness of the proposed scheme in simulating interesting test cases, such as water flow in heterogeneous soil and transient flow in variably saturated zones. Furthermore, the proposed semi-implicit scheme is applied to the coupled system of water flow and solute transport in unsaturated soils.

The rest of the paper is structured as follows. Section \ref{S2} introduces the coupled system of the Richards equation and the advection-dispersion equation. Moving on to Section \ref{S3}, the methodology and the studied temporal schemes are presented. In Section \ref{S4}, the investigation and comparison of these numerical schemes' convergence rates, accuracy, and robustness are conducted using exact and reference solutions, along with selected application tests. Concluding remarks are drawn in Section \ref{S5}.

%% file: Docts/ARichards_equation.tex
The Richards equation \cite{Richards1931} can be expressed in three distinct forms: the mixed form (\ref{equ21}), the $\Psi$-based form (\ref{equ22}), and the $\theta$-based form (\ref{equ23}). These forms are presented as follows  \cite{Celia1990}:\\
\begin{equation}
   \dfrac{\partial \theta}{\partial t}
= \nabla \cdot [ K_sk_r\nabla( \Psi  + z) ]+ s(\Psi),
   \label{equ21}
\end{equation}
\begin{equation}
 C(\Psi) \dfrac{\partial \Psi}{\partial t}
= \nabla \cdot [ K_sk_r\nabla( \Psi  + z) ]+ s(\Psi),
\label{equ22}
\end{equation}
\begin{equation} 
\dfrac{\partial \theta}{\partial t}
= \nabla \cdot [ D(\theta)\nabla\theta ] +  \dfrac{\partial K}{\partial z} + s(\theta).
\label{equ23}
\end{equation}
 In these equations, $\theta$  is the volumetric water content [$L^{3}/L^{3}$], $\Psi$ is the  pressure head [$L$], 
 $K_s(\bm{x})$ is the saturated hydraulic conductivity of soil [$L/T$], $k_r$ is the water relative permeability [$-$], $K = K_sk_r$ is  the unsaturated hydraulic conductivity [$L/T$], $z$ is the vertical coordinate positive upward [$L$], \(\bm{x}=(x,y,z)\) is the Cartesian coordinates of the domain, $t$ is the time [$T$], $s$ can be considered as  a sink or source term [$L^3/L^{3}/T$],
 \(C(\Psi)= \frac{\partial \theta}{\partial \psi}\) is the specific moisture capacity function [$1/L$] and \(D(\theta)= K\frac{\partial \Psi}{\partial \theta} \) is the unsaturated diffusivity [$L^2/T$].\\
 The flux of water $\textbf{q}$ in the soil is expressed as follows:
 $$\textbf{q}= -K\nabla( \Psi  + z).$$

%% file: Docts/Ph.RHC.tex
The water content of unsaturated soils is generally described in terms of the saturation $S$ [$-$]:
$$ 
S=\dfrac{\theta - \theta_r}{\theta_s - \theta_r}, 
$$
where $\theta_s$ and $\theta_r$ are the saturated water content and residual water content, respectively. Various empirical models have been developed to describe the saturation of unsaturated soils. These models provide valuable tools for understanding and predicting water flow behavior in unsaturated conditions. Among the commonly used models are the Gardner model \cite{GARDNER1958}, the Brooks-Corey model \cite{BrooksCorey1966}, the Haverkamp model \cite{havercamp1977} and the van Genuchten model (VG) \cite{VanGenuchten1980}. For the purpose of this study, the Gardner and VG models will be used. These models have been widely adopted in various studies \cite{Keita2021, Boujoudar2021, Boujoudar2023}. The  saturation can be written as follows for Gardner and van Genuchten models:\\
\begin{equation}
\begin{aligned}
    \textbf{Gardner model:} \quad S &= \exp(\alpha_v \Psi), \\
    \textbf{VG model:} \quad S &= \left[\dfrac{1}{1+(\alpha_v|\Psi|)^{n_{ v}}}\right]^{m_{ v}},
\end{aligned}
\label{equS}
\end{equation}
where $\alpha_v$ is the parameter related to the average pore size [$1/L$],  $n_{v}$  and $m_{ v}$ are dimensionless parameters related to pore size distribution.
For the water relative permeability, the Gardner and VG models \cite{Mualem1976}, are expressed as follows:\\
\begin{equation}
\begin{aligned}
    \textbf{Gardner model:} \quad k_r &= \exp(\alpha_v \Psi), \\
    \textbf{VG model:} \quad      k_r&=\dfrac{\left[1-(\alpha_v |\Psi|)^{n_{v}-1} (1+(\alpha_v |\Psi |)^{n_{v}})^{-m_{v}}\right]^2}{\left[1+\left(\alpha_v |\Psi |\right)^{n_{v}}\right]^{m_{v}/2}}.
\end{aligned}
\label{equK}
\end{equation}
 The VG's water relative permeability in \eqref{equK} is obtained according to Mualem's theory \cite{Mualem1976} where  $\quad m_{v} = 1 - \frac{ 1}{ n_{v}} \; \text{and} \; n_{v}>1.$
The pressure head ($\Psi \leq 0$) is describe using the Leverett $J$-function \cite{Levrett1941, Keita2021, Boujoudar2023}:
\begin{equation}
    \Psi=\psi_{c}J(S),
    \label{equ24}
\end{equation}
where $\psi_{c}$ is the capillary rise function ($\psi_{c} \geq 0)$) and $J$ is the  Leverett $J$-function ($J \leq 0$).  

%% file: Docts/ADE.tex
To describe the transport of non-reactive solutes in unsaturated soils, we usually use the mass conservation advection-dispersion equation. By ignoring the adsorption and chemical reactions, this equation, without sink or source terms, is expressed as follows:
\begin{equation}
   \dfrac{\partial \theta  c}{\partial t} - \nabla\cdot[ \theta \textbf{D}\nabla c- c\textbf{q}] = 0,
   \label{ADE}
\end{equation}
in which, $c$ represents the concentration rate in the liquid phase [$M/L^{3}$], $\textbf{D}$ is the dispersion tensor [$L^{2}/T$] and $\textbf{q}$ is the volumetric flux density [$L/T$]. The first term on the right side of \eqref{ADE} is the solute flux due to dispersion, the second term is the solute flux due to convection with flowing water. 
According to Bear \cite{bear2013dynamics} the components of the  dispersion tensor, \textbf{D}, in two-dimensional space (e.g. considering the  $x$ and $z$ directions) are defined as
\begin{equation}
\begin{aligned}
    &D_{xx}=\lambda_T |\textbf{v}| + \left(\lambda_L-\lambda_T\right)\dfrac{v_x^2}{|\textbf{v}|} + \tau(\theta)\lambda_m,\\
    &D_{xz}= \left(\lambda_L-\lambda_T\right)\dfrac{v_xv_z}{|\textbf{v}|},\\
    &D_{zz}=\lambda_T |\textbf{v}| + \left(\lambda_L-\lambda_T\right)\dfrac{v_z^2}{|\textbf{v}|} + \tau (\theta)\lambda_m,
\end{aligned}
\end{equation}
where $\lambda_L$ [$L$] and  $\lambda_T$ [$L$] are the longitudinal and transverse
dispersivities, $\textbf{v} :=\textbf{q}/\theta$ is the velocity vector, $\tau\left(\theta\right)$ is the tortuosity factor which depends on the water content, and $\lambda_m$ is the molecular diffusion coefficient in free water [$L^{2}/T$]. The tortuosity factor in the liquid phase is evaluated using the relationship of Millington and Quirk \cite{millington1961permeability}:
\begin{equation}
    \tau(\theta)=\dfrac{\theta^{7/3}}{\theta_s^2}.
\end{equation}

%% file: Docts/weak_form.tex
Let a computational domain $\Omega \subset \mathbb{R}^d$ ($d=2$) be an open-bounded domain with a sufficiently smooth boundary $\Gamma$, $T>0$ a given final computational time, $L^2(\Omega)$ represents the space of square-integrable real-valued functions defined on $\Omega$, $H^1(\Omega)$ denote its subspace that includes functions with first-order derivatives also in $L^2(\Omega)$. The mixed form of Richard's equation  without a sink or source term  can be expressed in terms of saturation as follows: 
\begin{equation}
   \phi \dfrac{\partial S}{\partial t}
- \nabla\cdot[ K\nabla( \Psi  + z) ] = 0.
   \label{equ31}
\end{equation}
Multiplying \eqref{equ31} and \eqref{ADE} by a test function $v\in  H^1(\Omega)$ and integrating over the domain, we obtain:
\begin{equation}
\left\{
    \begin{array}{ll}
\displaystyle\int_{\Omega}\phi\dfrac{\partial S}{\partial t}v\;\mathrm{d}\bm{x} + \int_{\Omega}[K\nabla( \Psi  + z)]\cdot\nabla v\;\mathrm{d}\bm{x} +\int_{\Gamma}q_w v\; \mathrm{d} \Gamma=0,
\\
\\
\displaystyle\int_{\Omega}\dfrac{\partial (\theta c)}{\partial t}v\;\mathrm{d}\bm{x} + \int_{\Omega} \left[\theta \textbf{D} \nabla c - c \textbf{q}\right]\cdot\nabla v\;\mathrm{d}\bm{x}    +\int_{\Gamma}q_c v\; \mathrm{d} \Gamma=0.
    \end{array}
\right.
\label{equ32}
\end{equation}
Here, $ \phi:= \theta_s-\theta_r$ and the two variables $q_w$ and $q_c$ $\in \mathbb{R}$ represent the water and solute fluxes across the boundary $\Gamma$, respectively, and are expressed as follows:
\begin{equation}
\left\{
    \begin{array}{ll}
q_w = \textbf{q} \cdot \mathbf{n}   \quad  &\text{on} \quad \Gamma,
\\
\\
q_c  \;= -\left(\theta \textbf{D} \nabla c - c \textbf{q}\right) \cdot \mathbf{n}  \;\quad  &\text{on} \quad \Gamma,

    \end{array}
\right.
\end{equation}
where $\mathbf{n}$ is the outward normalized vector on the boundary of the domain. A weak formulation of the system \eqref{equ31}-\eqref{ADE} is stated as: Determine $\left(S(\cdot, t),c(\cdot, t)\right) \in H^1(\Omega)\times H^1(\Omega)$ such that the system \eqref{equ32} is satisfied for every $v \in H^1(\Omega)$ and a.e $t\in (0, T]$.\\
In this study, we will employ the standard finite element method to showcase the significance of the proposed time discretization approach.
Let $\mathcal{T}_h$ being an unstructured partitioning of the domain $\Omega$ into disjoint triangles $\kappa$, so that
\begin{equation}
    \overline{\Omega}
 = \bigcup_{\kappa \in \mathcal{T}_h} \kappa. \label{eq:domain_partition}
\end{equation}
The element mesh size $h$ is defined as
\begin{equation}
    h := \max_{\kappa \in \mathcal{T}_h} h_\kappa, \label{eq:mesh_size}
\end{equation}
where $h_\kappa:= \text{diam}(\kappa)$ and let $P_k$ represent the polynomial space with $k$ denoting the spatial polynomial degree associated with ${\mathcal{T}_h}$, which will be taken equal to 1 in this study. The standard finite element space is defined as:
$$ \mathcal{
V}_h = \left\{ v_h \in H^1(\Omega) : v_h|_{K_j} \in P_k(K_j), \forall K_j \in \mathcal{T}_h \right\}. $$
By using the linear Galerkin finite elements defined above in space, the
semi-discrete variational formulation of \eqref{equ32} reads as:
\begin{equation}
\left\{
    \begin{array}{lll}
    \text{ {\fontfamily{qtm}\selectfont{ For a.e. }}} t \in (0,T] \text{ {\fontfamily{qtm}\selectfont{ and }}} v_h \in \mathcal{V}_h,
\text{ {\fontfamily{qtm}\selectfont{ find }}} (S_h, c_h) \in \mathcal{V}_h\times \mathcal{V}_h \text{ {\fontfamily{qtm}\selectfont{ such that there holds }}}
\\[7pt]
\displaystyle\int_{\Omega}\phi_h\dfrac{\partial S_h}{\partial t}v_h\;\mathrm{d}\bm{x} + \int_{\Omega}[K_h\nabla( \Psi_h  + z)]\cdot\nabla v_h\;\mathrm{d}\bm{x} +\int_{\Gamma}q_w v_h \;\mathrm{d}\Gamma=0, 
\\[7pt]
\displaystyle\int_{\Omega}\dfrac{\partial (\theta_h c_h)}{\partial t}v_h\;\mathrm{d}\bm{x} + \int_{\Omega} [\theta_h \textbf{D}_h \nabla c_h - c_h \textbf{q}_h]\cdot\nabla v_h\;\mathrm{d}\bm{x}    +\int_{\Gamma}q_c v_h\; \mathrm{d} \Gamma=0,
    \end{array}
\right.
\label{equ3w}
\end{equation}
In this study, we will use both Dirichlet and Neumann boundary conditions. On some boundaries of the domain, Dirichlet boundary conditions are applied. Under these conditions, the solutions $S$ and $c$ of the governing equations are prescribed at the boundaries. Furthermore, the test function $v$ is constrained to be zero on these boundaries.

%% file: Docts/timediscrete.tex
Let  $\Delta t>0$ represents the time step and the closed interval $\left[0, T\right]$ discretized into a collection of $N > 0$ subintervals, so that,   \begin{align}
    t_n = n \Delta t, \quad n = 0, 1, \ldots, N, \quad \Delta t = \dfrac{T}{N}.
\end{align}
We define $u_h^n$ as the approximation to $u_h$ at time level $t_n$, where \(u_h\) represents the solution of the system \eqref{equ3w}. In this study, we propose a second-order time stepping method to discretize in time the system \eqref{equ3w}. This method includes two free parameters $\delta$ and $\mu \in [0, 1]$. For a given problem
\begin{equation}
\dfrac{\partial\varphi}{\partial t}=F(\varphi,t),
\label{model_equ}
\end{equation}
the proposed time discretization technique is obtained by centering the scheme in time about the time level $t_{n+\delta}$ as follows:
\begin{equation}
\dfrac{(\delta +\frac{1}{2}) \varphi^{n+1}-2\delta \varphi^{n}+(\delta -\frac{1}{2})\varphi^{n-1}}{\Delta t} = (\delta+\mu)F^{n+1}+(1-\delta-2\mu) F^n+\mu F^{n-1},
\label{Gene_equ}
\end{equation}
where \(\varphi^n\) and \(F^{n}\) denote the solution and the right-hand side term of equation \eqref{model_equ} calculated at time level $t_n$, respectively. The idea behind this technique of time discretization is to generalize the BDF2 scheme with more accurate and effective schemes. In connection with this, Beljadid et al. \cite{Beljadid2014} studied the schemes \eqref{Gene_equ} with $\delta = 1$ and proposed a modified BDF2 method along with a range of other methods. It was, specifically, used to tackle the linear part of some problems applicable to atmospheric models. In our analysis, we will specify some schemes resulting from the general one \eqref{Gene_equ}. By setting \(\delta=1\) and \(\mu=0\) , the scheme is BDF2 method. By setting \(\delta=1\) and \(\mu=1 \), the scheme will be referred to as a 2nd-order semi-implicit backward differentiation formulae scheme (SBDF2). For \(\delta=\frac{1}{2}\) and \(\mu=0\), the scheme corresponds to the 2nd-order Crank Nicholson method and will be denoted by CN2. As a summary of the above extracted methods, the time discretization of the Richards and solute transport equations using these proposed methods can be expressed with the weak formulation \eqref{equ3w} as follows: Given an appropriate approximation of the initial solution  $\left(S_h^0, c_h^0\right) \in \mathcal{V}_h\times \mathcal{V}_h$ and a proper initialization for $(S_h^1, c_h^1) \in \mathcal{V}_h\times \mathcal{V}_h$, find $(S_h^{n+1},c_h^{n+1}) \in \mathcal{V}_h\times \mathcal{V}_h$ such that:\\
\textbf{BDF2:} 
\begin{equation}
\left\{
\begin{array}{lll}
\displaystyle\int_{\Omega}\phi_h\left(\dfrac{3S_h^{n+1}-4S_h^{n}+S_h^{n-1}}{2\Delta t}\right) v_h\;\mathrm{d}\bm{x} 
+\int_{\Omega}\left[K\left( \Psi_h^{n+1} \right)\nabla\left( \Psi_h^{n+1}  + z\right)\right]\cdot\nabla v_h\;\mathrm{d}\bm{x}
+\int_{\Gamma}q_w v_h \;\mathrm{d} \Gamma=0,\\[5pt]
\displaystyle\int_{\Omega} \left(\dfrac{3\theta_h^{n+1} c_h^{n+1} -4\theta_h^{n} c_h^{n}+\theta_h^{n-1} c_h^{n-1}}{2\Delta t}\right)v_h\;\mathrm{d}\bm{x}+\displaystyle\int_{\Omega} \left[\theta_h^{n+1} \mathbf{D}_h^{n+1} \nabla c_h^{n+1}
- \textbf{q}_h^{n+1}c_h^{n+1}\right] \cdot \nabla v_h\;\mathrm{d}\bm{x} \\[5pt]
+\displaystyle\int_{\Gamma}q_c v_h\; \mathrm{d} \Gamma=0, \qquad \forall v_h \in \mathcal{V}_h.
\label{BDF2}
\end{array}
\right.
\end{equation}
\textbf{SBDF2:}\\
\begin{equation}
\left\{
\begin{array}{lll}
\begin{aligned}
&\int_{\Omega}\phi_h\left(\dfrac{3S_h^{n+1}-4S_h^{n}+S_h^{n-1}}{2\Delta t}\right) v_h\;\mathrm{d}\bm{x} 
 +2\int_{\Omega}\left[K( \Psi_h^{n+1} )\nabla( \Psi_h^{n+1}  + z)\right]\cdot\nabla v_h\; \mathrm{d}\bm{x} -2\int_{\Omega}[K\left( \Psi_h^{n} \right)\\[5pt]
 &\times\nabla\left( \Psi_h^{n}  + z\right)]\cdot\nabla v_h\; \mathrm{d}\bm{x}
+\int_{\Omega}\left[K \left( \Psi_h^{n-1} \right)\nabla\left( \Psi_h^{n-1} + z\right)\right]\cdot\nabla v_h\;\mathrm{d}\bm{x} +\int_{\Gamma}q_w v_h \;\mathrm{d} \Gamma=0,\\[5pt]
&\displaystyle\int_{\Omega} \left(\dfrac{3\theta_h^{n+1} c_h^{n+1} -4\theta_h^{n} c_h^{n}+\theta_h^{n-1} c_h^{n-1}}{2\Delta t}\right)v_h\;\mathrm{d}\bm{x}+2\displaystyle\int_{\Omega} \left[\theta_h^{n+1} \mathbf{D}_h^{n+1} \nabla c_h^{n+1}
- \textbf{q}_h^{n+1}c_h^{n+1}\right] \cdot \nabla v_h\;\mathrm{d}\bm{x}\\[5pt]
& -2\displaystyle\int_{\Omega} \left[\theta_h^{n} D_h^{n} \nabla c_h^{n} - \textbf{q}_h^{n}c_h^{n}\right] \cdot \nabla v_h\;\mathrm{d}\bm{x}
+\displaystyle\int_{\Omega} \left[\theta_h^{n-1} \mathbf{D}_h^{n-1} \nabla c_h^{n-1}
- \textbf{q}_h^{n-1}c_h^{n-1}\right] \cdot \nabla v_h\;\mathrm{d}\bm{x}
\\[5pt]
&+\displaystyle\int_{\Gamma}q_c v_h \;\mathrm{d} \Gamma=0, \qquad \forall v_h \in \mathcal{V}_h. \qquad \qquad \qquad \qquad \qquad \qquad \qquad \qquad \qquad \qquad \qquad \qquad \qquad \qquad \qquad \qquad \qquad \qquad
\end{aligned}
\label{SBDF2}
\end{array}
\right.
\end{equation}
\textbf{CN2:}\\
\begin{equation}
\left\{
\begin{array}{lll}
\begin{aligned}
&\int_{\Omega}\phi_h\left(\dfrac{S_h^{n+1}-S_h^{n}}{\Delta t}\right)v_h\;\mathrm{d}\bm{x} 
+\dfrac{1}{2}\left(\int_{\Omega}\left[K\left( \Psi_h^{n+1} \right)\nabla\left( \Psi_h^{n+1}  + z\right)\right]\cdot\nabla v_h\;\mathrm{d}\bm{x}\right)
 \\[5pt]
 &+\dfrac{1}{2}\left(\int_{\Omega}\left[K\left( \Psi_h^{n} \right)\times\nabla\left( \Psi_h^{n}  + z\right)\right]\cdot\nabla v_h\;\mathrm{d}\bm{x}\right)
+\int_{\Gamma} q_w v_h \;\mathrm{d} \Gamma=0,\\[5pt]
 &\displaystyle\int_{\Omega} \left(\dfrac{\theta_h^{n+1} c_h^{n+1} -\theta_h^{n} c_h^{n}}{\Delta t}\right)v_h\;\mathrm{d}\bm{x}+\frac{1}{2}\left(\displaystyle\int_{\Omega} \left[\theta_h^{n+1} \mathbf{D}_h^{n+1} \nabla c_h^{n+1}
- \textbf{q}_h^{n+1}c_h^{n+1}\right] \cdot \nabla v_h\;\mathrm{d}\bm{x}\right)\\[5pt]
&+\frac{1}{2}\left(\displaystyle\int_{\Omega} \left[\theta_h^{n} D_h^{n} \nabla c_h^{n}
- \textbf{q}_h^{n}c_h^{n}\right] \cdot \nabla v_h\;\mathrm{d}\bm{x}\right)+\displaystyle\int_{\Gamma}q_c v_h \;\mathrm{d} \Gamma=0, \qquad \forall v_h \in \mathcal{V}_h. \qquad \qquad \qquad \qquad \qquad \qquad
 
\end{aligned}
\end{array}
\right.
\label{ACN2}
\end{equation}
Another method that can be derived from \eqref{Gene_equ} is the Leapfrog scheme, which is centered at the time level \(t_n\). This method is obtained by setting \(\delta=0\) and \(\mu=0\). However, the Leapfrog scheme is explicit and can pose stability challenges, necessitating the introduction of a stabilization term. This stabilization term uses a second-order approximation for \(\Psi_h^{n}\) and \(c_h^{n}\) in the Richards and advection-dispersion equations, respectively. The approximations are given by:
\begin{equation}
\begin{cases}
    \Psi_h^n &\simeq \nu \Psi_h^{n+1} + (1-2\nu)\Psi_h^n + \nu \Psi_h^{n-1}, \\
    c_h^n &\simeq \nu c_h^{n+1} + (1-2\nu)c_h^n + \nu c_h^{n-1}, 
\end{cases}
\label{approximationSILF2}
\end{equation}
where \(0<\nu\leq 1\) is a stabilization parameter. The resulting semi-implicit method is referred to as SILF2 scheme.
The weak formulation \eqref{equ3w} using SILF2 scheme reads as follows: Given an initial approximation \((S_h^0, c_h^0) \in \mathcal{V}_h \times \mathcal{V}_h\) and an appropriate initialization for \((S_h^1, c_h^1) \in \mathcal{V}_h \times \mathcal{V}_h\), find \((S_h^{n+1}, c_h^{n+1}) \in \mathcal{V}_h \times \mathcal{V}_h\) such that:\\
\textbf{SILF2:}\\
\begin{equation}
\left\{
\begin{array}{lll}
\begin{aligned}
&\int_{\Omega}\phi_h \dot{S}^n_h\left(\dfrac{ \Psi_h^{n+1}-\Psi_h^{n-1}}{2\Delta t}\right)v_h \;\mathrm{d}\bm{x} 
 +\int_{\Omega}\left[K( \Psi_h^{n} )\nabla\left(\Psi_h^{n} +\nu\left(\Psi_h^{n+1} -2\Psi_h^{n}+ \Psi_h^{n-1}\right)+z\right)\right]\cdot \nabla v_h \;\mathrm{d}\bm{x}\\[5pt]
 &+\int_{\Gamma}q_w v_h \;\mathrm{d} \Gamma=0,
 \\[5pt]
 & \theta^{n+1} = \left(\theta_s - \theta_r\right)\left(1+(\alpha_v|\Psi^{n+1}|)^{n_{ v}}\right)^{-m_{ v}} + \theta_r,\\[5pt]
&\displaystyle\int_{\Omega} \left(\dfrac{\theta_h^{n+1} c_h^{n+1} - \theta_h^{n-1} c_h^{n-1}}{2\Delta t}\right)v_h\;\mathrm{d}\bm{x}+\displaystyle\int_{\Omega} \Big[\theta_h^n \mathbf{D}_h^n \nabla\left(c_h^n+ \nu[c_h^{n+1}-2c_h^{n}+c_h^{n-1}]\right)\\[5pt]
&
- \textbf{q}_h^n \cdot \left(c_h^n +\nu[c_h^{n+1}-2c_h^{n}+c_h^{n-1}]\right) \Big] \cdot \nabla v_h \;\mathrm{d}\bm{x}+\displaystyle\int_{\Gamma}q_c v_h\; \mathrm{d} \Gamma,\qquad \forall v_h \in \mathcal{V}_h,
 \end{aligned}
 \end{array}
\right.
\label{SILf}
\end{equation}
where $\mathbf{q}^n = -K\left(\Psi_h^n, \bm{x}\right)\nabla\left(\Psi_h^n + z\right)$ and $\dot{S}_h^n:=\frac{\partial S^n_h}{\partial \Psi_h}$. Note that for simplicity, we write \(K(\Psi_h)\) for the unsaturated hydraulic conductivity which depends also on space (\(K(\Psi_h, \bm{x})\)). The primary idea behind introducing the SILF2 method centers on:
\begin{itemize}
    \item Using the \(\Psi\)-based form where the nonlinear terms \(C(\Psi_h^n)\) and \(K(\Psi_h^n)\) in \eqref{equ22} are known at time level \(t_n\), thereby simplifying the system to a linear state with \(\Psi_h^{n+1}\) as the primary variable of interest. This process involves discretizing both the left- and right-hand sides of equation \eqref{equ22} at the specific time level \(t_n\).
    \item This discretization specifically results in assigning values of \(\delta = 0\) and \(\mu = 0\) in \eqref{Gene_equ},  consequently transforming the right-hand side of equation \eqref{equ22} into an explicit formulation. To address potential instability issues that may arise from the approach's explicit formulation, we incorporate a stabilization term governed by the adaptable free parameter \(\nu\). This addition is crucial for enhancing the method's robustness and ensuring computational stability.
\end{itemize}
Semi-implicit methods are widely used in the literature to solve time-dependent partial differential equations including linear and nonlinear operators \cite{ascher1995implicit,  JIANG2021607, KEITA2021107588, Keita2021}. The linear part is integrated implicitly while the nonlinear one is integrated explicitly \cite{ascher1995implicit}. However, since the Richards equation is a nonlinear PDE, we opted to integrate the nonlinear operator with a combination of implicit and explicit terms, leading to the scheme \eqref{Gene_equ}. 
Concerning the error analysis in time of these schemes and as we will show in \ref{appA:details}, the method \eqref{Gene_equ} converges to \eqref{model_equ} with a second-order truncation error for all \(\delta\) and \(\mu\). Therefore, BDF2, SBDF2, and CN2 methods exhibit a second-order accuracy in time. For the SILF2 method, the approximations \eqref{approximationSILF2} reveal a second-order accuracy for all \(\nu\) (see \ref{appA:details}). As a result, the SILF2 method achieves second-order accuracy in time.
\begin{remark}
    In cases where the medium is saturated (\( S = 1\)), the term \(\dot{S}\) in the SILF2 method will vanish. This leads to the disappearance of the first term in the discretized Richards equation in \eqref{SILf}.
    While the unknown variable \(\Psi^{n+1}\) is eliminated in this term, it remains in the second term of the discretized Richards equation.  Numerical investigations have been conducted in the following sections in the presence of saturated zones, and the applied SILF2 approach has been found to yield stable results.
    \label{remark3.2}
\end{remark}
\begin{remark}
     To initiate the process, we consider the initial conditions $S^0_h=\mathbf{P}_hS_0 \in V_h$ and $c^0_h=\mathbf{P}_hc_0 \in V_h$, where $\mathbf{P}_h :H^1_0(\Omega)\rightarrow V_h$ denotes the standard projection operator, \(S_0\) and \(c_0\) are the initial conditions added to the governing equation of unsaturated flow and solute transport, respectively. The second step of the solution $\left(\Psi_h^1,c_h^1\right)$ is calculated employing the implicit Backward Euler method (BDF1):
     \begin{equation}
     \left\{
\begin{array}{lll}
\begin{aligned}
&\int_{\Omega}\phi_h \left(\dfrac{S_h^{1}-S_h^{0}}{\Delta t}\right)v_h\;\mathrm{d}\bm{x} 
 +\int_{\Omega}[K( \Psi_h^{1} )\nabla( \Psi_h^{1}  + z)]\cdot\nabla v_h\;\mathrm{d}\bm{x} +
\int_{\Gamma}q_w v_h \;\mathrm{d} \Gamma=0,\\[5pt]
&\displaystyle\int_{\Omega} \left(\dfrac{\theta_h^{1} c_h^{1} -\theta_h^{0} c_h^{0}}{\Delta t}\right)v_h\;\mathrm{d}\bm{x}+\displaystyle\int_{\Omega} \left[\theta_h^{1} \mathbf{D}_h^{1} \nabla c_h^{1}
-\textbf{q}_h^{1}c_h^{1}\right] \cdot \nabla v_h\;\mathrm{d}\bm{x}
+\displaystyle\int_{\Gamma}q_c v_h \;\mathrm{d} \Gamma=0, \quad \forall v_h \in \mathcal{V}_h.
\label{equ39}
 \end{aligned}
 \end{array}
\right.
\end{equation} 
This approach necessitates the application of an iterative technique for solving the first equation of the system. We use the modified Picard method as explained in the next section to linearize the equation. It is important to note that this linearization technique will be utilized in the first step of the SILF2 method, but it will be uniformly employed at all steps for the other methods which are nonlinear.
\end{remark}

%% file: Docts/linearisation.tex
\begin{figure}
\centering
\includegraphics[scale=0.5]{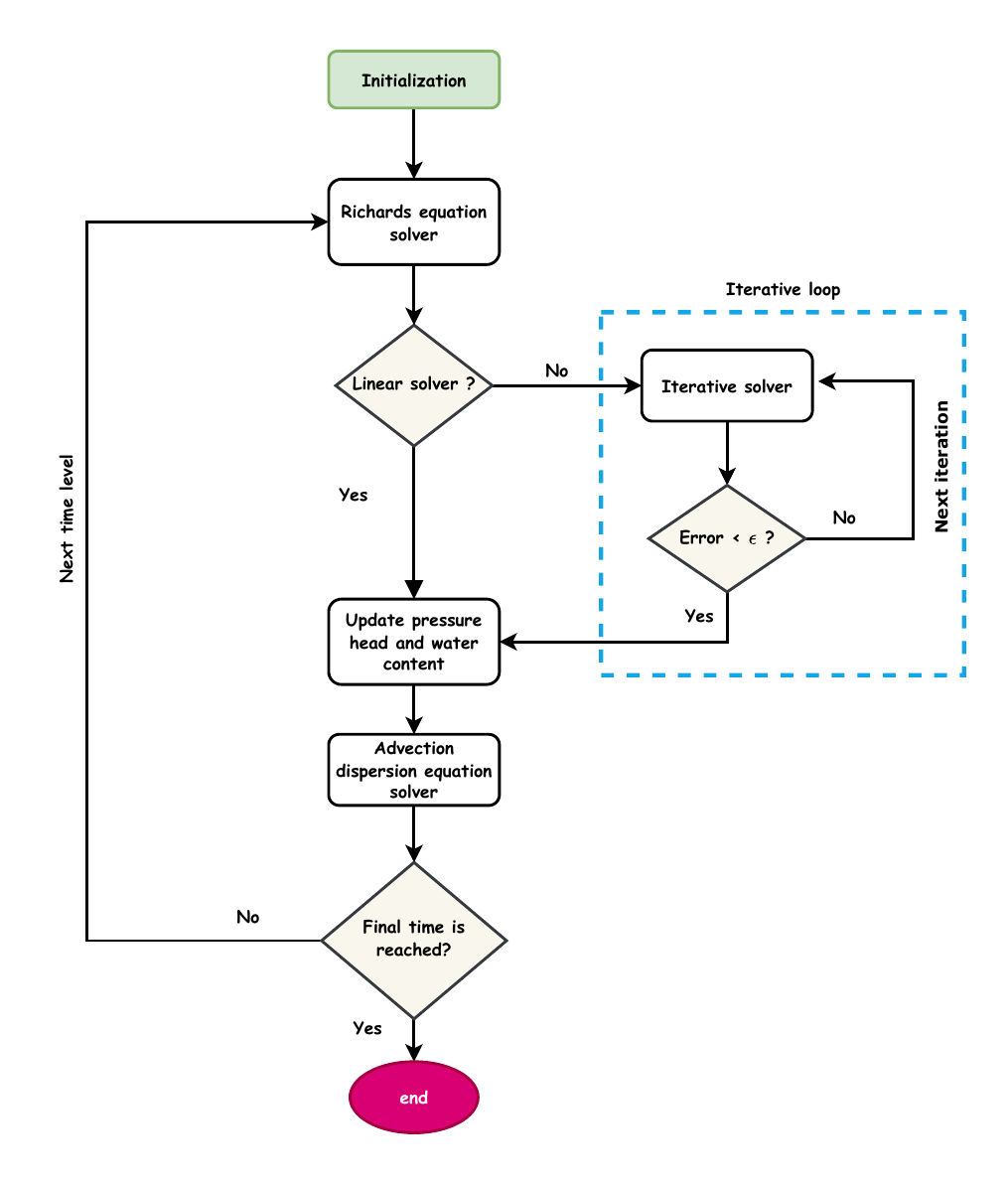} 
\caption{Flowchart diagram for methodology implementation.}
\label{Fig:flowshare}
\end{figure}

As we mentioned before, Richards' equation exhibits strong nonlinearity due to the presence of the hydraulic conductivity term and the specific moisture
capacity term. Therefore, to address the complexity of Richards' equation with fully implicit schemes, linearization techniques are essential. Several iterative methods have been employed to linearize the equation, such as the Newton method, Picard method, and L-scheme method, among others \cite{albuja2021family, Celia1990, radu2016, Lott2012, Mitra2019}. These techniques play a crucial role in making the solution to the equation more tractable and efficient. Several studies have extensively examined these linearization methods by comparing them with each other and adapting strategies to deal with convergence issues. Understanding the performance and limitations of each method has enabled researchers to propose adaptive strategies that can effectively overcome convergence challenges.
For instance, Mitra et al. \cite{Mitra2019} provide a modified L-scheme that combines the L-scheme and the modified Picard methods. List and Radu \cite{radu2016} used the L-scheme/Newton method, which is a combination of the L-scheme method and Newton one, and the results show good performance in terms of robustness and time computation.
Within the framework of this study, we employ the methodology of the modified Picard approach to linearize the Richards equation in its mixed form, which is widely used in the literature \cite{Celia1990, Huang1996, radu2016}. The modified Picard method utilizes Taylor series expansion to approximate the current iteration level of saturation $S^{n+1, m+1}$ \cite{Celia1990}, so that,
\begin{equation}
    S^{n+1, m+1}=S^{n+1, m}+\dot{S}^{n+1, m}\left(\Psi^{n+1, m+1}-\Psi^{n+1, m}\right)+O\left((\Psi^{n+1, m+1}-\Psi^{n+1, m})^2\right).
    \label{equ39}
\end{equation}
The Newton method exhibits quadratic convergence, while the modified Picard method demonstrates linear convergence. Although Newton's advantageous convergence properties, it requires an initial guess that is so close to the solution to achieve this convergence. The modified Picard method is easy to implement and demands less memory compared to Newton’s method. 
In each iteration, we start the process by initializing from the previous time step. In other words, \(\Psi^{n+1,0}\) is set equal to \(\Psi^{n}\).  To stop the iteration process, many studies have utilized different stopping criterion expressions \cite{Celia1990, clement2021, Huang1996,radu2016}. Within the scope of this study, we will use the following criterion:
\begin{equation}
    || \Psi^{n+1,m+1}-\Psi^{n+1,m} ||_{L^2} \leq \epsilon,
\end{equation}
in which \(0 < \epsilon \ll 1\)
 is the tolerance parameter and \(|| \cdot ||_{L^2}\) represents the L\(^2\)-norm, which is defined as
 \begin{equation}
     || u ||_{L^2}=\left(\int_{\Omega} |u|^2\, \mathrm{d}\bm{x}\right)^{\frac{1}{2}},
     \label{L2Error}
 \end{equation}
for a given function \(u \in L^2(\Omega)\).
\begin{remark}
Celia et al. \cite{Celia1990} studied the finite element methods for solving the Richards equation and illustrated that they require the use of mass lumping on the time derivative to effectively prevent nonoscillatory solutions. Different forms of mass lumping have been developed for standard and mixed finite element methods to model infiltration in porous media \cite{Belfort2009, Celia1990, Younes2006}. In our study, the mass lumping technique used in \cite{Celia1990} is applied for the system of Richards equation and the solute transport equation.
\end{remark}
\begin{remark}
In this study, the coupled system of unsaturated flow and solute transport in porous media is considered a \textquotedblleft one-way coupling\textquotedblright  system, where the dynamic of water influences the transport process. In this context, the Richards equation for water flow is solved first to obtain the flow field, and then the advection-dispersion equation is solved for the solute transport. The two equations are connected through this sequential process, representing one-way influence of the flow on solute transport as shown in Figure \ref{Fig:flowshare}.
\end{remark}

%% file: Docts/numerical_tests.tex
\begin{figure}
\begin{center}
\includegraphics[width=\linewidth]{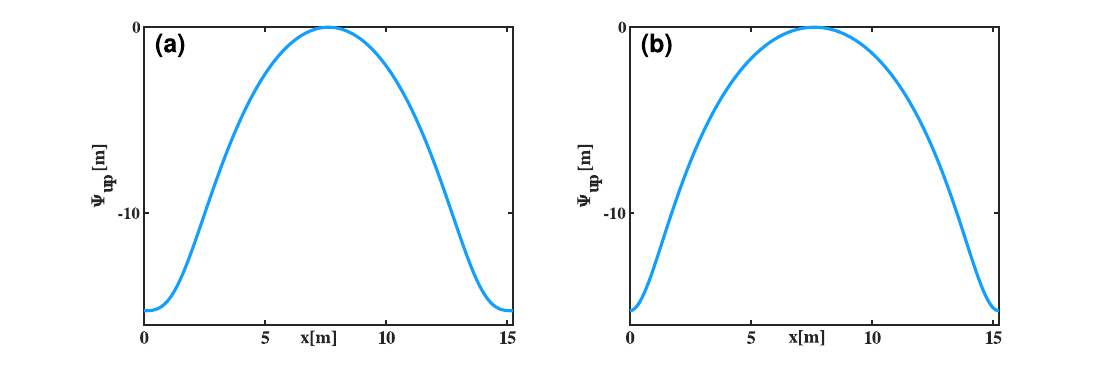} 
\caption{The boundary condition enforced on the top of the soil. \textbf{(a)} and \textbf{(b)} refer to Test 1 and Test 2, respectively.}
\label{Fig:2dBC}
\end{center}
\end{figure}
In this part, the proposed methods will be evaluated using the Gardner model for the capillary pressure. To achieve this, we will leverage analytical solutions from reputable sources in the literature \cite{Boujoudar2021, Keita2021, tarcy2011}. The capillary rise function and the Leverett $J$-function are expressed as follows:
\begin{equation}
    \psi_c=\dfrac{1}{\alpha_v}, \quad \quad J(S)=\log(S).
\end{equation}

\setlength{\arrayrulewidth}{1pt}
\noindent
\begin{table}[ht]
\centering

\begin{tabularx}{\textwidth}{>{\hsize=2\hsize} XX >{\hsize=0.2\hsize} XX}
\hline
\textbf{Variable} & \textbf{Symbol} & \textbf{Value}  \\
\hline
Length of the square domain & $L$ [$m$]& $15.24$ \\
Hydraulic conductivity & $K_s$ [$m/day$]& $0.10$ \\
Saturated water content & $\theta_s$ \;\,[$m^3/m^3$] & $0.450$  \\
Residual water content & $\theta_r$ \;\,[$m^3/m^3$] & $0.150$  \\
Gardner fitting coefficient  & $\alpha_v$ \;\;\,[$1/m$]& $0.164$ \\
Initial condition  & $\Psi_d$ \,[$m$]& $-15.24$  \\
Tolerance criterion  & $\epsilon$ \;\;\,\,\,[$-$]& $10^{-6}$  \\
Series terms number  & $nt$ \;\;[$-$]& $200$  \\
\hline
\end{tabularx}
\caption{Parameters and material properties used for 2D infiltration tests.}
\label{tab:3}
\end{table}
The two numerical tests are performed using the Green and Ampt problem \cite{green_ampt} with different boundary conditions.  The initial soil is set extremely dry (\(\Psi=-15.24\, m\)).
In Figure \ref{Fig:2dGA}\textbf{a} and \ref{Fig:2dGA}\textbf{b}, we depict the precise top boundary condition that defines the soil infiltration procedure. For simplicity, we will refer to the two tests as Test 1 for the first test problem and Test 2 for the second test problem. As illustrated in Figure \ref{Fig:2dGA}, the boundary conditions of Test 1 and Test 2 are expressed in the following formulations:
\begin{enumerate}
    \item[] Test 1: Dirichlet boundary conditions are enforced on the four sides:
        \begin{equation}
\left\{
    \begin{array}{ll}
\begin{aligned}
\Psi(x,z,t=0)& =\Psi(x,z=0,t)=  \Psi(x=0,z,t)=\Psi(x=L,z,t)=\Psi_d,\\
\Psi(x,z=L,t)&=\Psi_{up}(x)=\frac{1}{\alpha_v}\log\left(\zeta +\left(1-\zeta \right)\left[\frac{3}{4}\sin\left(\frac{\pi x}{L}\right)-\frac{1}{4}\sin\left(\frac{3\pi x}{L}\right)\right]\right).    
\end{aligned}
    \end{array}
    \right.
\end{equation}
\item[] Test 2: No flux boundary conditions are enforced on the lateral sides of the domain, while Dirichlet boundary conditions are prescribed at the upper and bottom boundaries:
        \begin{equation}
\left\{
    \begin{array}{ll}
\begin{aligned}
\Psi(x,z,t=0)& =\Psi(x=0,z,t)=\Psi_d,\\
\Psi(x,z=L,t)&=\Psi_{up}(x)=\frac{1}{\alpha_v}\log\left[\zeta +\frac{\left(1-\zeta \right)}{2}(1-\cos\left(\frac{2\pi x}{L}\right)\right],   
\end{aligned}
    \end{array}
    \right.
\end{equation}
\end{enumerate}
where \(\Psi_{up}(x)\) is a function that depends on $x$ and is described as the top boundary condition of the two tests. The two-dimensional analytical solutions of the Richards equation related to these tests are expressed as follows: \cite{Boujoudar2021, Keita2021, tarcy2011}:
      \begin{equation}
    \Psi(x,z,t)=\psi_c \log(\zeta + \Psi_0(x,z,t)),
\end{equation}
where \(\Psi_0(x,z,t)\) is  given for each test by the following expressions:\\
Test 1: 
\begin{equation}
\left\{
    \begin{array}{ll}
\begin{aligned}
    & \Psi_0(x,z,t)=(1-\zeta)\exp(\frac{\alpha_v}{2}(L-z))\bigg(\dfrac{3}{4}\sin(\dfrac{\pi x}{L})\bigg[\dfrac{\sinh(\beta_1 z)}{\sinh(\beta_1 L)}
    +\dfrac{2}{Ld}\sum_{p=1}^{\infty}(-1)^{^p} \dfrac{\lambda_p}{\nu_{1p}}\sin(\lambda_p z)\\
    &\times \exp(-\nu_{1p} t)\bigg]
    -\dfrac{1}{4}\sin(\dfrac{3\pi x}{L})\bigg[\dfrac{\sinh(\beta_3 z)}{\sinh(\beta_3 L)}
    +\dfrac{2}{Ld}\sum_{p=1}^{\infty}(-1)^{^p} \dfrac{\lambda_p}{\nu_{3p}}\sin(\lambda_p z)\exp(-\nu_{3p} t)\bigg]\bigg).
    \end{aligned}
    \end{array}
    \right.
\end{equation}
Test 2: 
\begin{equation}
\begin{aligned}
    & \Psi_0(x,z,t)=(1-\zeta)\sin\left(\frac{\pi x}{L}\right)\exp(\frac{\alpha_v}{2}(L-z))\bigg[\dfrac{\sinh(\beta_1 z)}{\sinh(\beta_1 L)}
    +\dfrac{2}{Ld}\sum_{p=1}^{\infty}(-1)^{^p} \dfrac{\lambda_p}{\nu_{1p}}\sin(\lambda_p z)
    \times \exp(-\nu_{1p} t)\bigg]
    \end{aligned}
\end{equation}
with
\begin{equation}
    \zeta=\exp(\alpha_v\Psi_d),\; d=\dfrac{\alpha_v\phi}{K_s}, \; \lambda_p=\dfrac{\pi p}{L},\; \beta_i=\Big(\frac{\alpha_v^2}{4}+(\frac{i\pi}{L})^2\Big)^{\frac{1}{2}},\; \text{and} \; \nu_{ip}=\dfrac{\beta_i^2 +\lambda_p^2}{d} \quad .
\end{equation}
We numerically solved the Richards equation in the square domain $\Omega=[0, L]\times [0, L]$ shown in Figure \ref{Fig:2dGA} using the studied numerical methods. The parameters and soil proprieties used in the numerical tests are depicted in Table \ref{tab:3} \cite{Boujoudar2021, Keita2021}.

\begin{figure}
\centering
\includegraphics[width=\linewidth, height = 18cm]{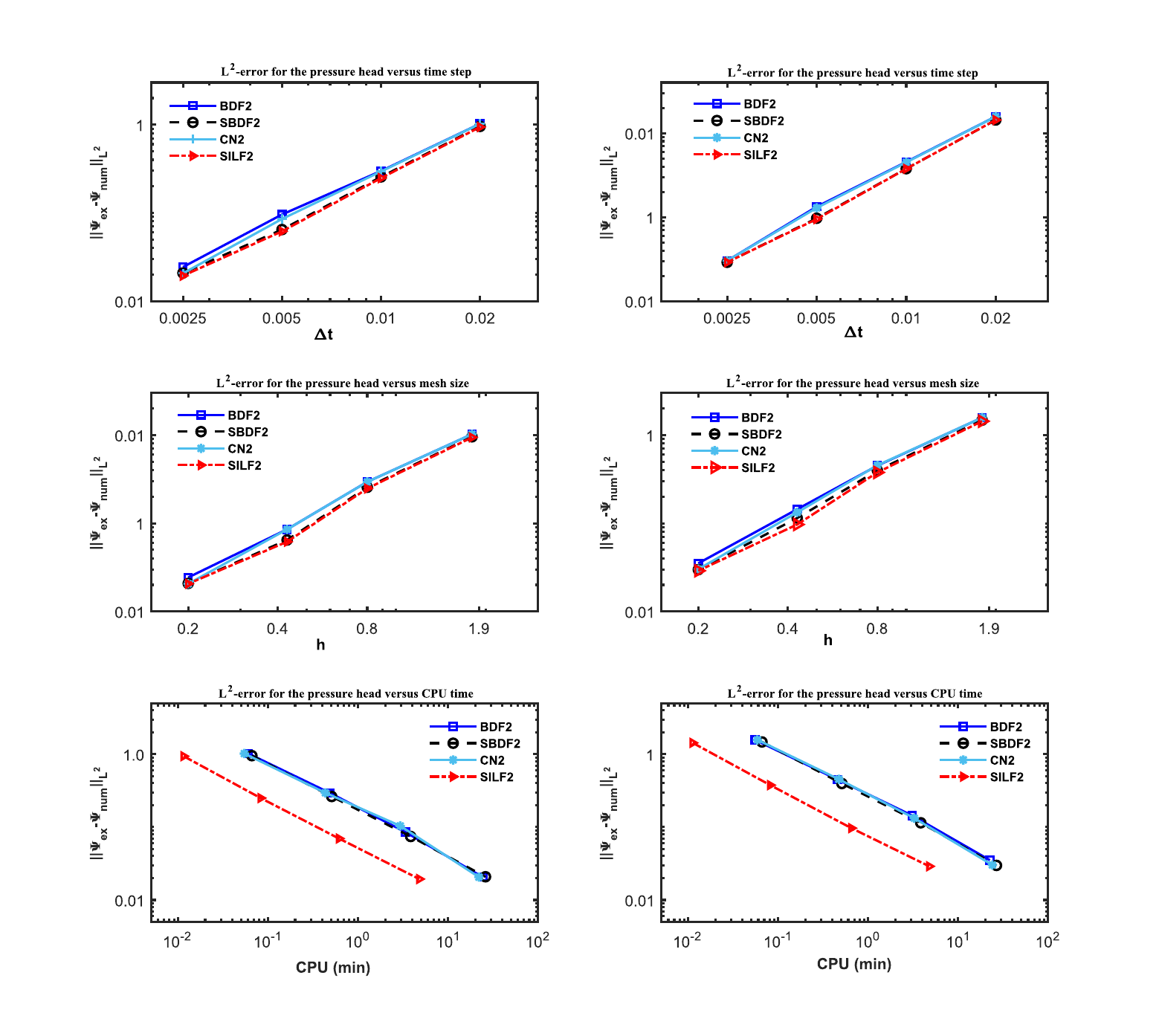} 
\caption{$L^2$ error on log-scaled of the pressure head for Test 1 (left column) and Test 2 (right column) as a function of time step (top), mesh size (middle), and CPU time (bottom) for all proposed methods.}
\label{Errors}
\end{figure}
\begin{figure}
\begin{center}
\includegraphics[width=\linewidth]{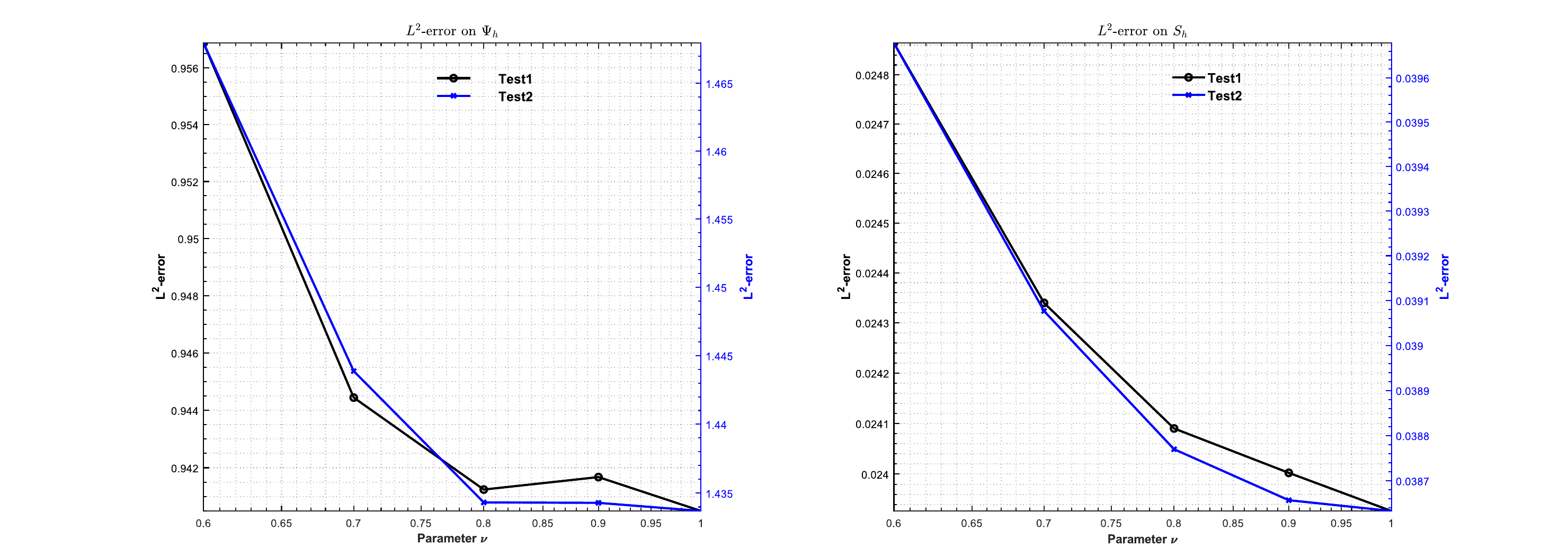}
\caption[]{$L^2$-error for Test 1 and Test 2 on log-log scale using the SILF2 method with different values of the free parameter \(\nu\).}
\label{Fig:parameternu}
\end{center}
\end{figure}
\begin{figure}
\begin{center}
\includegraphics[width=0.45\linewidth]{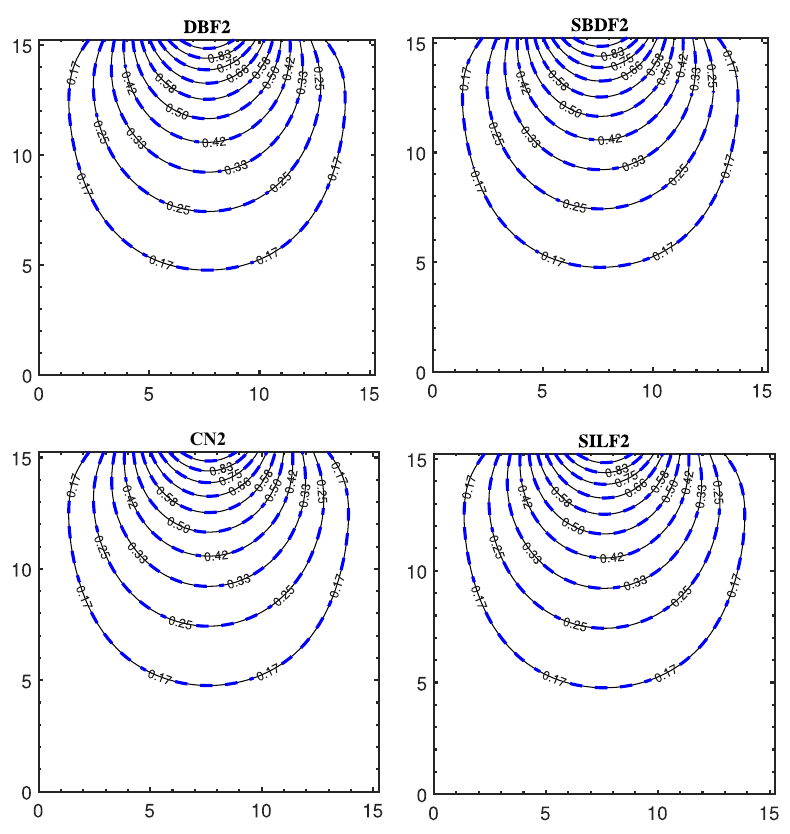}
\hspace{1.2cm}
\includegraphics[width=0.45\linewidth]{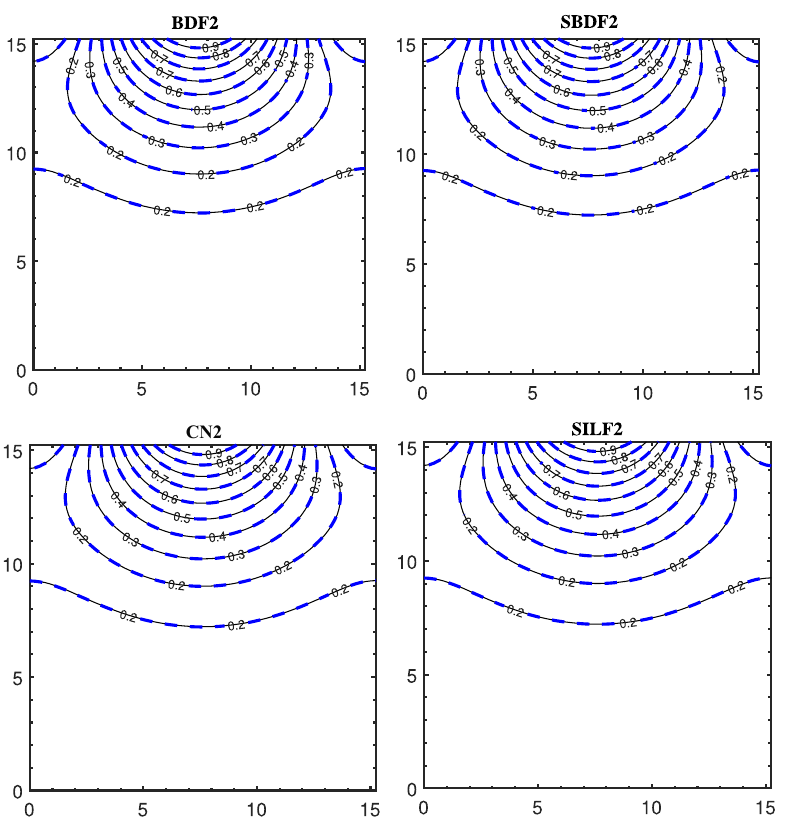}
\caption[]{Numerical results for Test 1 (left) and Test 2 (right). Comparison between numerical solution (\tikz[baseline=-0.5ex]\draw[blue, dashed, thick] (0,0) -- (1,0);) using BDF2, SBDF2, CN2 and SILF2 methods and exact solution  (\tikz[baseline=-0.5ex]\draw[black, thick] (0,0) -- (1,0);) for the saturation.}
\label{Fig:Contoure_plot}
\end{center}
\end{figure}

In the numerical analysis, we opt to study the $L^2$-error versus time step, versus mesh size, and versus CPU time for a final time $T=5$ days with a varying time step of $\Delta t =0.02,\, 0.01,\, 0.005,\, 0.0025$ and mesh size of $h=1.9,\, 0.8,\, 0.4,\, 0.2$, where the subdivisions $12\times12,\, 25\times25,\, 50\times50$ and $100\times100$ along the boundaries of the domain are used respectively. A suitable tolerance criterion was selected to ensure the convergence of the iterative technique and the accuracy of the methods. However, this can lead to some observed computational overhead regarding CPU time. Figure \ref{Errors} depicts the $L^2$-error values of the pressure head obtained with varying \(\Delta t,\, h\) and CPU for both two tests. We analyze the SILF2 method using different values of the free parameter \(\nu\) to ensure both stability and accuracy. As shown in Figure \ref{Fig:parameternu}, the optimal value of the free parameter is \(\nu = 1\). It is important to note that our analysis focuses on the parameter \(\nu\) within the range of 0.6 to 1. Values of \(\nu\) between 0 and 0.5 are excluded from the analysis due to convergence issues. 
The SILF2 method with this optimal free parameter is the most accurate scheme, followed by the SBDF2 method. Regarding the computational time, the SILF2 method leads to accurate results with less computational cost thanks to its linear nature starting from the second step. In terms of accuracy, all methods are of order $2$ demonstrating a satisfactory agreement with the theoretical analysis presented in \ref{appA:details}. Table \ref{tab:comparison_case1} offers a detailed analysis of $L^2$-error, rate of convergence, and CUP time of the BDF2, SBDF2, CN2, and SILF2 methods. The comparison between the spatial distribution of the saturation for the numerical and analytical solutions at the final time $T=5$ days for both tests is illustrated in Figure \ref{Fig:Contoure_plot}. For the numerical solution, we utilized a mesh size of $h= 0.4$ ($50\times50$) and a time step of $\Delta t=0.005$. The contour plots showcase good agreement between the numerical and analytical solutions. 

\begin{table}[htbp]
    \centering
    \caption{$L^2$-error, order of convergence and CPU time as a function of mesh size and time step of the studied methods for the two numerical tests.}
    \small 
    \setlength{\tabcolsep}{12pt} 
    \renewcommand{\arraystretch}{1} 
    \begin{tabular}{ccccccccccc}
    \toprule
   Method & $h$ & $\Delta t$ & \multicolumn{3}{c}{Test 1} & \multicolumn{3}{c}{Test 2} \\
    \cmidrule(lr){4-6} \cmidrule(lr){7-9} \cmidrule(lr){9-11}
     &  &   & $L^2$-error & Order & CPU (s) & $L^2$-error   & Order & CPU (s) \\
    \midrule
  BDF2&1.9&2E-2&1.02326 & &3.62 &1.57566 & &3.32\\
      &0.8&1E-2&0.2982  & &29.13 &0.44962 & &27.36\\
      &0.4&5E-3&0.095769  & &201.89 &0.144199 & &184.95\\
      &0.2&2.5E-3&0.0243305 &$1.97$ &1460.54 &0.03518 &$2.03$ &1351.67\\
      \\
SBDF2&1.9&2E-2&0.956127 & &3.96 &1.4757 & &3.96\\
      &0.8&1E-2&0.262932 & &30.55 &0.396407 & &30.63\\
      &0.4&5E-3&0.0742976 & &230.88 &0.114618 & &230.54\\
      &0.2&2.5E-3&0.0208178 &$1.84$ &1564.93 &0.0297928 &$1.94$ &1600.08\\
      \\
   CN2&1.9&2E-2&1.02183 & &3.27 &1.57495 & &3.52\\
      &0.8&1E-2&0.295176 & &26 &0.453597 & &28.32\\
      &0.4&5E-3&0.103785 & &176.77 &0.132803 & &194.94\\
      &0.2&2.5E-3&0.0206168 &$2.09$ &1322.78 &0.0301976 &$2.1$ &1440.91\\
      \\
SILF2&1.9&2E-2&0.940499 & &0.69 &1.43371 & &0.67\\
      &0.8&1E-2&0.250411 & &4.96 &0.376912 & &4.89\\
      &0.4&5E-3&0.0696979 & &36.86 &0.0968422 & &39.02\\
      &0.2&2.5E-3&0.0193712 &$1.85$ &287.72 &0.0289613 &$1.86$ &286.7\\      

    \bottomrule
    \end{tabular}
    \label{tab:comparison_case1}
\end{table}

\begin{figure}
\begin{center}
\includegraphics[width=\linewidth]{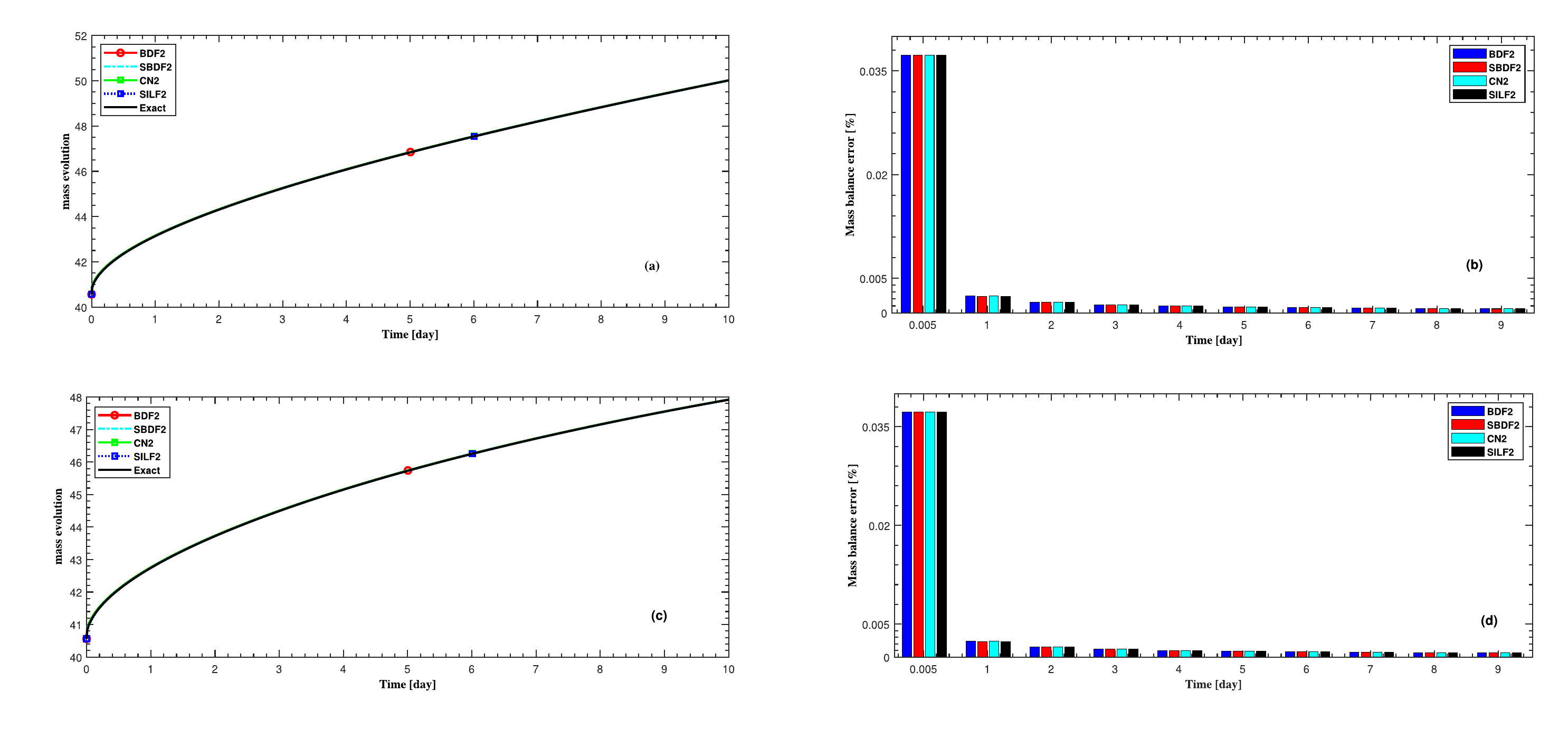}
\caption{Total mass analysis for Test 1 (left) and Test 2 (right). (\textbf{a}) and (\textbf{b}) represent the total mass evolution over time of the numerical and exact solutions. (\textbf{c}) and (\textbf{d}) represent the relative mass balance errors of the numerical methods at discrete times $t = 7.2$ min, $t= 1$ day, $t= 2$ days, $\cdots$, $t= 9$ days. }
\label{Fig:MB}
\end{center}
\end{figure}  

In the following, we compare the total mass evolution for the numerical and  exact solutions using the following terms:
\begin{equation}
    \int_{\Omega} \theta(\bm{x},t)\, \mathrm{d}\bm{x} \quad \text{and} \quad \int_{\Omega} \theta_{ex}(\bm{x},t)\, \mathrm{d}\bm{x},
\end{equation}
where \(\theta\) and \(\theta_{ex}\) represent the numerical and exact solutions of the moisture content, respectively. For all the studied methods, we computed the relative mass balance error (MBE) at each time level using the following formulae \cite{Celia1990, ROMANO1998315}:
\begin{equation}
    \text{MBE}(t)=\Bigg| 1- \dfrac{\text{MB}_{num}(t)}{\text{MB}_{ex}(t)}\Bigg|\times 100,
    \label{eqError}
\end{equation} 
where  \(\text{MB}_{num}(t)=\int_{\Omega} \left[\theta(\bm{x}, t)-\theta_0(\bm{x})\right] \, \mathrm{d}\bm{x} \) and \(\text{MB}_{ex}(t)=\int_{\Omega} \left[\theta_{ex}(\bm{x}, t)-\theta_0(\bm{x})\right] \, \mathrm{d}\bm{x} \)  represent the total mass change that occurs in the domain \( \Omega\) within the time period for the numerical and exact solutions, and \(\theta_0(\bm{x}) := \theta(\bm{x}, t=0)\) is the initial moisture content. We consider the mesh size $h=0.8$ where a total of 5000 triangles and 2601 vertices are used, the time step $ \Delta t = 0.005$ and the final time $T=10$ days. Figure \ref{Fig:MB}\textbf{a} and \ref{Fig:MB}\textbf{c} depict the resulting mass evolution for the numerical and exact solutions for Test 1 and Test 2, respectively. The results demonstrate good agreement between the numerical mass evolution and the exact one for both tests. Figure \ref{Fig:MB}\textbf{b} and \ref{Fig:MB}\textbf{d} illustrate the corresponding relative mass balance error for Test 1 and Test 2, respectively.
The results indicate a mass error at the beginning of the evolution, which diminishes rapidly as the period of evolution extends.

%% file: Docts/Heter.tex
\begin{figure}
\begin{center}
\includegraphics[scale=0.8]{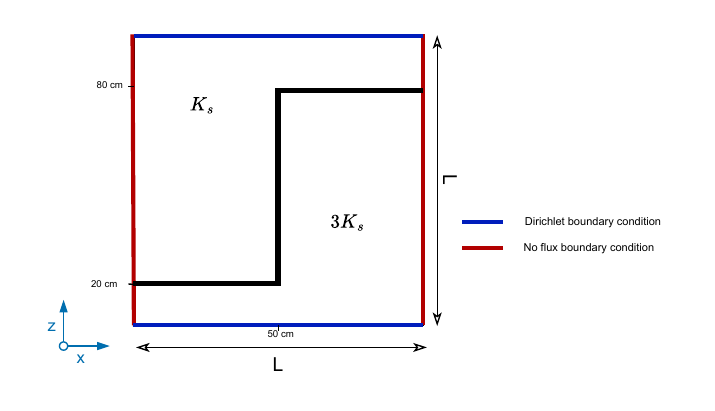}
\caption{Schematic illustration of the L-shape domain separated into two sub-regions.}
\label{Fig:lshape}
\end{center}
\end{figure}
\begin{figure}
\begin{center}
\includegraphics[width=\linewidth]{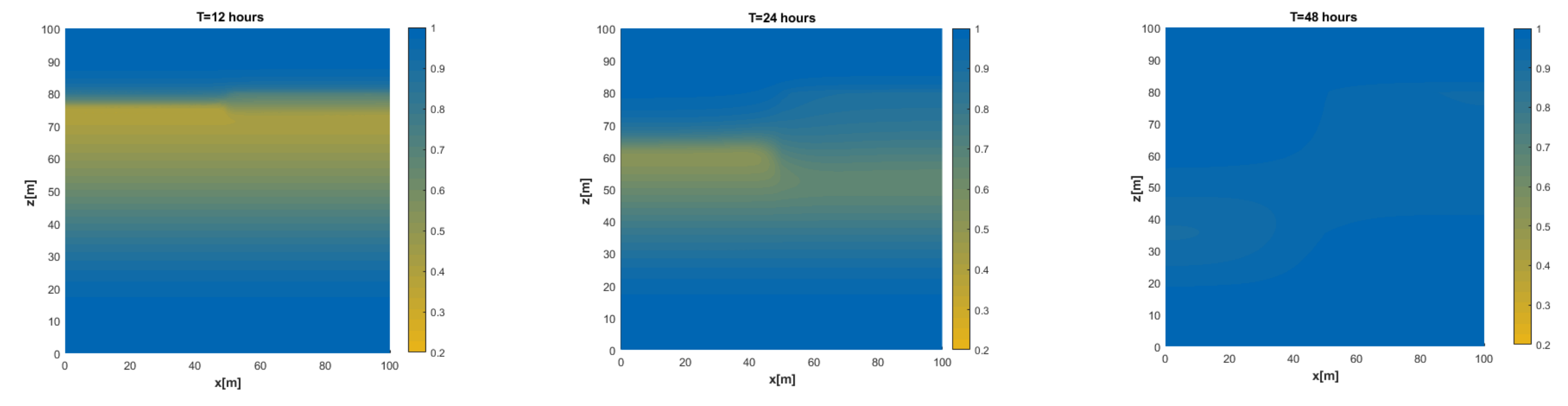}
\caption{Time evolution of the saturation using SILF2 scheme for the heterogeneous test. }
\label{Fig:Lshape}
\end{center}
\end{figure}

The previous part related to 2D homogeneous soils confirmed that the SILF2 method is accurate and it requires a small CPU time. In this part, we aim to test the efficiency and performance of this method in two-dimensional heterogeneous soil. The VG model is used to describe the relation between the saturation, the pressure head and the hydraulic conductivity. The capillary rise and the Leverett $J$-function are given by: 
\begin{equation}
    \psi_c=\dfrac{1}{\alpha_v}, \quad \quad J(S)=-\big[S^{-\frac{1}{m_{v}}} - 1\big]^{\frac{1}{n_{v}}}.
    \label{equJ}
\end{equation} 
A square computational domain with a length of \(100\;cm\) is separated into two sub-regions by an L-shape form as illustrated in Figure \ref{Fig:lshape}. We consider the initial condition $\Psi(x, z, t=0) = -z$. Homogeneous Dirichlet boundary conditions are used at the top and bottom boundaries, while homogeneous Neumann boundary conditions are imposed at the lateral sides.
\setlength{\arrayrulewidth}{1pt}
\noindent
\begin{table}[ht]
\centering

\begin{tabularx}{\textwidth}{>{\hsize=2\hsize} XX >{\hsize=0.4\hsize} XX}
\hline
\textbf{Material propriety} & \textbf{Symbol} & \textbf{value}  \\
\hline
Saturated hydraulic conductivity  & $K_s$[$cm/s$] & $6.944E-5$ \\
Saturated water content & $\theta_s$ [$cm^3/cm^3$] \;\, & $0.5$  \\
Residual water content & $\theta_r$ [$cm^3/cm^3$] & $0.12$  \\
VG fitting coefficient   & \(\alpha_v \)  [$1/cm$] & $0.02$ \\
VG parameter & $n_v$[$-$] & $3$ \\
\hline
\end{tabularx}
\caption{Material properties used for the heterogeneous test \cite{Keita2021}.}
\label{tab:6}
\end{table}\\
This test was carried out in  \cite{BARON2017104}, where the two-layered soils of the domain are considered to have the same features except the saturated hydraulic conductivity. The first soil has a saturated hydraulic conductivity $K_s$, while the second has a saturated hydraulic conductivity equal to $3K_s$ as illustrated in Figure \ref{Fig:lshape}. The material properties related to this test are given in Table \ref{tab:6}.  The domain is discretized with a nonuniform triangular mesh with a mesh size of \(h=1\, cm\) where a total of 93234 elements and 47045 vertices are used. The time step \(\Delta t = \frac{1}{180} \) hour is chosen, and the simulation time is $T = 48$ hours. Figure \ref{Fig:Lshape} demonstrates the evolution of the saturation in time simulated using the SILF2 scheme and the results agree well with those of the previous studies for this test \cite{BARON2017104, Boujoudar2023,  Keita2021}. Regarding the convergence analysis, a reference solution is computed by discretizing the domain with a total of $1254412$ triangles and $628780 $ vertices and with a time step \(\Delta t = \frac{1}{360} \) hour. Here, we use the relative $L^2$-error, which is given by 
\begin{equation}
    \frac{\| u_{\text{ref},h} - u_h \|_{L^2}}{\| u_{\text{ref},h} \|_{L^2}},
\end{equation}
where $u_{\text{ref},h}$ and $u_h$ represent the reference and numerical solutions , respectively.
Table \ref{tab:Error_Lshape} displays the relative \(L^2\)-errors for both the pressure head and saturation as well as the CPU time in function of the varying time step. An order of convergence around $1.63$  is observed for the pressure head.

\begin{table}[htbp]
    \centering
    \caption{Relative $L^2$-error and CPU time as a function of time step using the SILF2 method for the heterogeneous soil test.}
    \setlength{\tabcolsep}{12pt} 
    \renewcommand{\arraystretch}{1} 
    \begin{tabular}{cccc}
    \toprule
    $\Delta t$(s) & Relative $L^2$-error on $\Psi_h$ & Relative $L^2$-error on $S_h$ & CPU(s) \\
    \midrule
   $160$ & 0.183642 & 0.083599 & 5041.71 \\
   $80$ & 0.13213 & 0.0502052 & 5888.65 \\
   $40$ & 0.0831433 & 0.0223597 & 10035.9 \\
   $20$ & 0.0524077 & 0.00721635 & 12056 \\
    \bottomrule
    \end{tabular}
    \label{tab:Error_Lshape}
\end{table}

%% file: Docts/Vaclintest.tex
\begin{figure}[htbp]
        \includegraphics[width=0.48\linewidth]{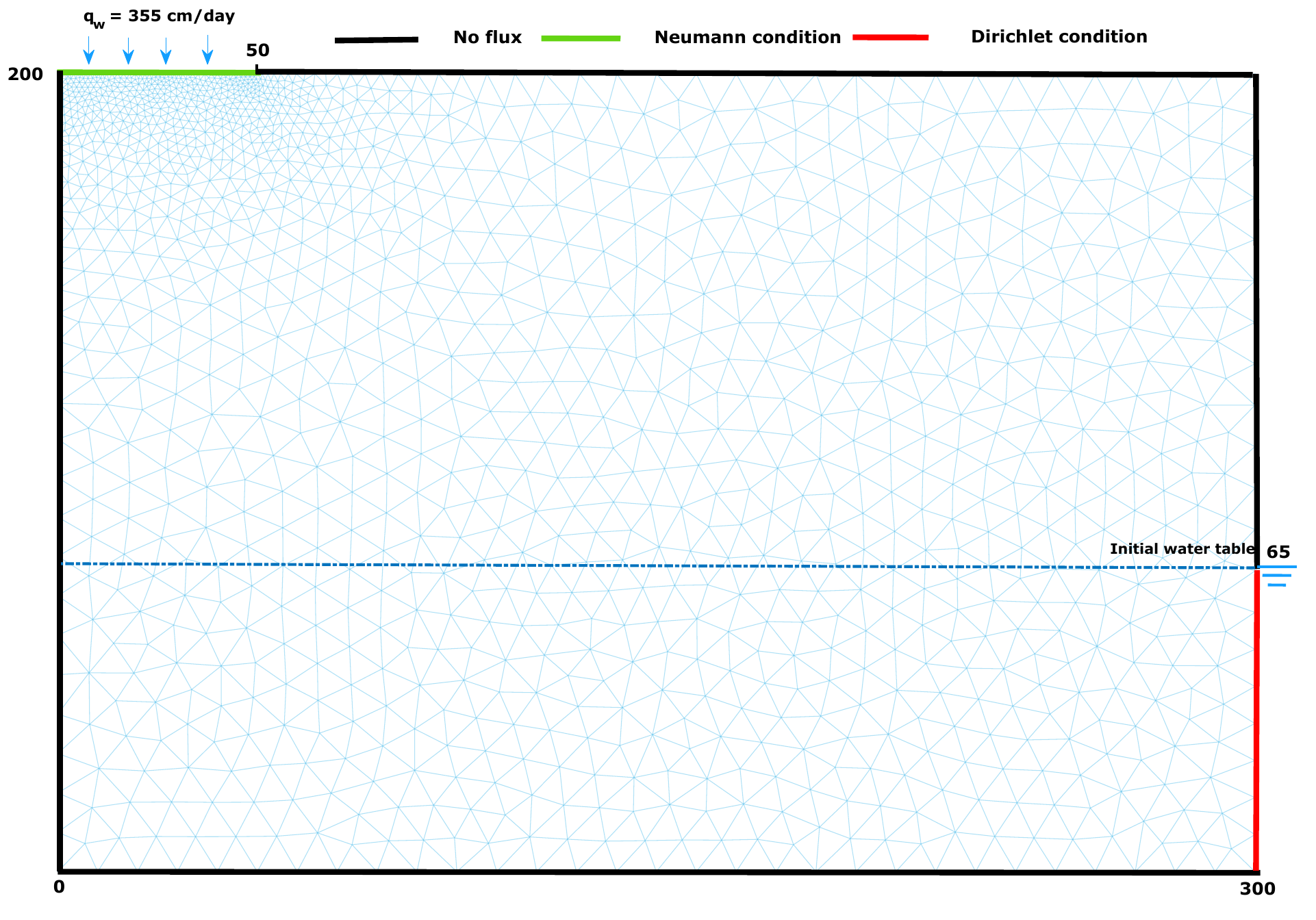}
        \includegraphics[width=0.49\linewidth]{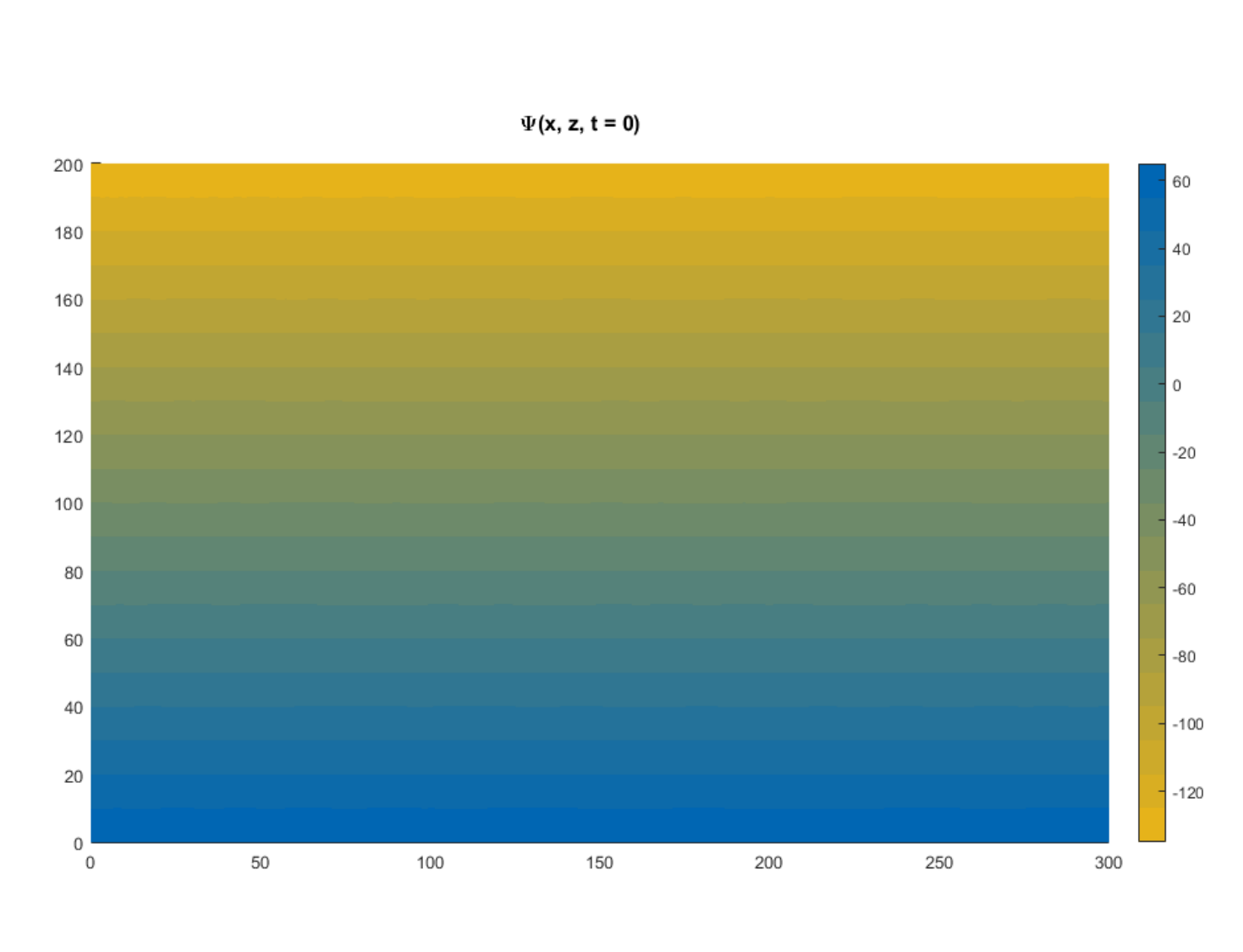}
    \caption{Schematic illustration of the right half of the domain (left) and initial pressure head distribution (right) for the water table recharge test case.}
    \label{fig:Vauclin_configuration}
\end{figure}
\setlength{\arrayrulewidth}{1pt}
\noindent
\begin{table}[ht]
\centering
\begin{tabularx}{\textwidth}{>{\hsize=2\hsize} XX >{\hsize=0.4\hsize} XX}
\hline
\textbf{Material propriety} & \textbf{Symbol} & \textbf{value}  \\
\hline
Saturated hydraulic conductivity  & $K_s$ [$cm /day$] & $840$ \\
Saturated water content & $\theta_s$ [$cm^3/cm^3$] \;\, & $0.3$  \\
Residual water content & $\theta_r$ [$cm^3/cm^3$] & $0.01$  \\
VG fitting coefficient   & \(\alpha_v \)  [$1/cm$] & $0.033$ \\
VG parameter   & \(n_v\)  [$-$] & $4$ \\
\hline
\end{tabularx}
\caption{Material proprieties used for water table recharge test-case \cite{CLEMENT199471, FAHS20091122}. }
\label{tab:Vauclin}
\end{table}

\begin{figure}
\centering
\includegraphics[width=0.6\linewidth]{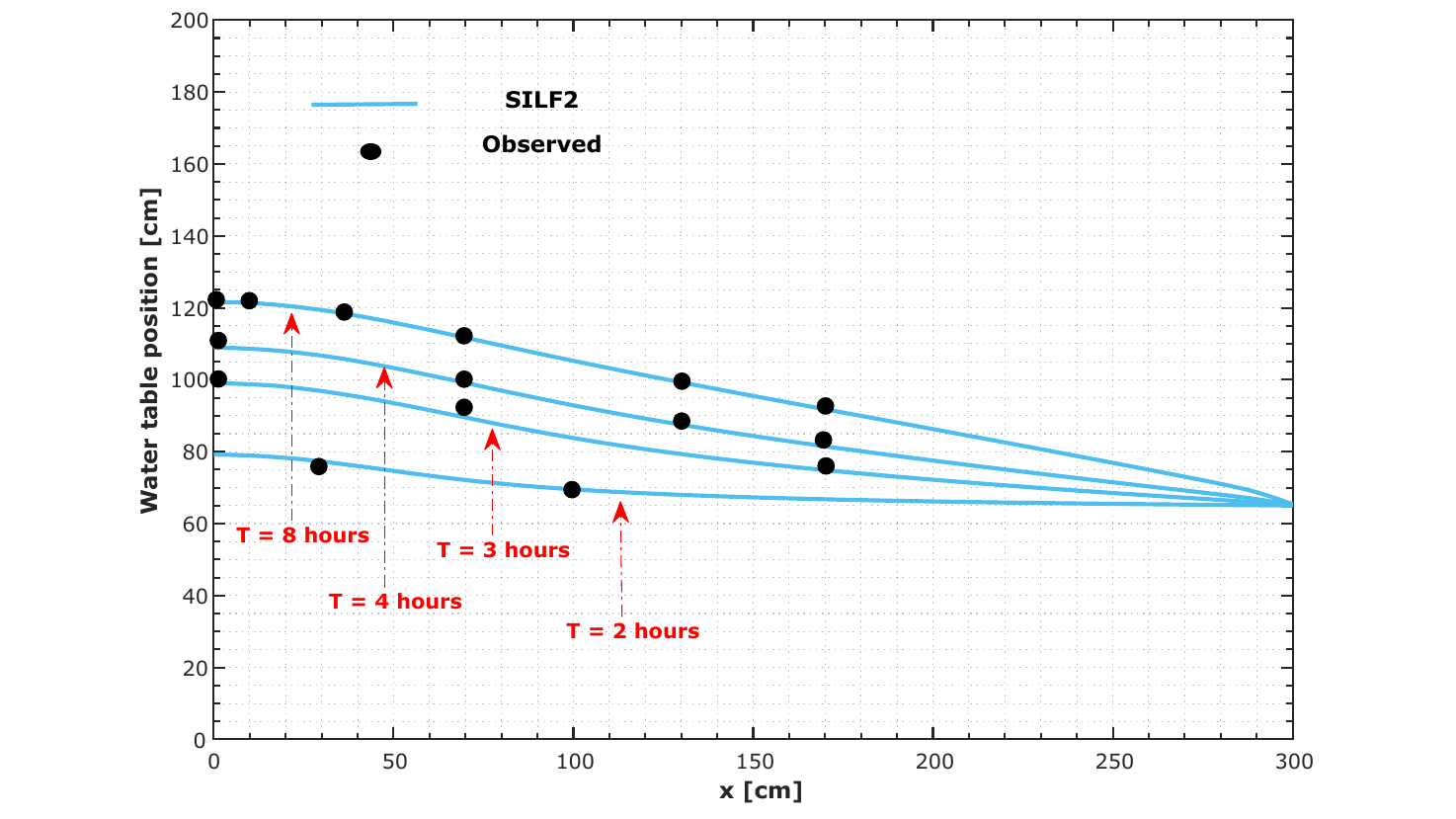}
\caption{Comparison of simulation and experimental data \cite{vauclin1979experimental} for groundwater levels at various time intervals.}
\label{fig:water_table}
\end{figure}

\begin{figure}
\centering
\includegraphics[width=\linewidth]{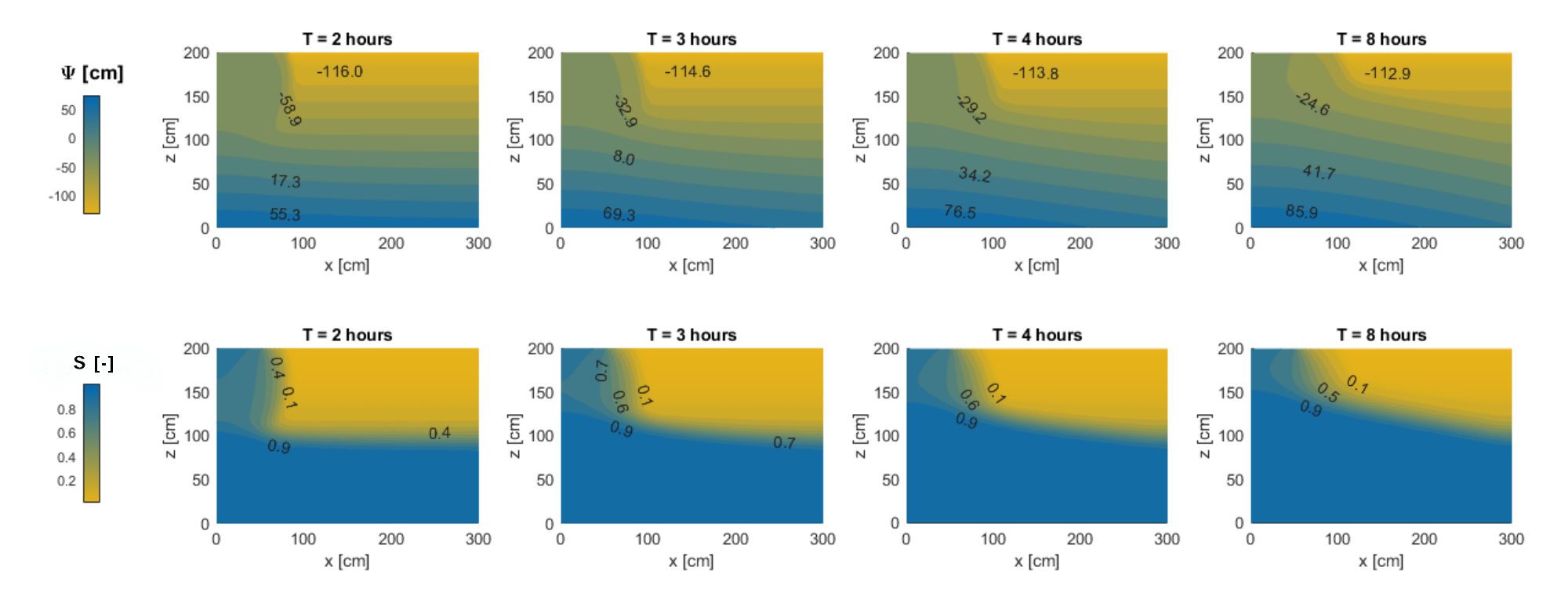}
\caption{Time evolution of the pressure head (top) and saturation (bottom)  for the water table recharge test case using the SILF2 scheme.}
\label{fig:evolution}
\end{figure}

In this numerical test, we evaluate the performance of the proposed SILF2 method for simulating variably saturated soil, based on the experiments carried out in \cite{vauclin1979experimental}.  A rectangular computational domain of 600 $cm$ in width and 200 $cm$ in height with the initial water table at 65 $cm$ above the base is used for the experimental setup. At the soil surface, a recharge rate of 355 $cm/day$ is applied over a central width of 100 $cm$. Due to the symmetry of the problem, only the right half of the domain is simulated. The boundaries included no-flow conditions along the symmetry axis, the upper boundary (except in the recharge zone), and the lower boundary. The right boundary had a fixed hydrostatic pressure of 65 $cm$ below the initial water table, with a no-flow condition above it. Initially, the soil system was in hydrostatic equilibrium with the water table. Figure \ref{fig:Vauclin_configuration} shows the test case configuration and initial condition. Hydraulic properties are expressed using the VG relations in equations \eqref{equS} and \eqref{equK}. The material properties, sourced from previous studies \cite{CLEMENT199471, FAHS20091122}, are detailed in Table \ref{tab:Vauclin}. A triangular mesh with 2643 elements and 1406 vertices is used with mesh refinement carried out near the recharge zone as shown in Figure \ref{fig:water_table}. A time step $\Delta t = 0.036$ seconds is used, and the evolution of the levels of water table is computed at 2, 3, 4, and 8 hours during the simulations. Figure \ref{fig:water_table} compares the experimental and simulated water table levels. The evolution of the water table levels is accurately depicted, closely matching the experimental profiles reported in \cite{vauclin1979experimental}. Furthermore, Figure \ref{fig:evolution} shows the simulated solutions of pressure head and saturation at time $t = 2$ hours, $t = 3$ hours, $t = 4$ hours, $t = 6$ hours, and $t = 8$ hours. The results show the capability of SILF2 method to simulate water flow in variably saturated soils.


%% file: Docts/SolutePENNS.tex
\begin{figure}
\begin{center}
\includegraphics[scale = 1]{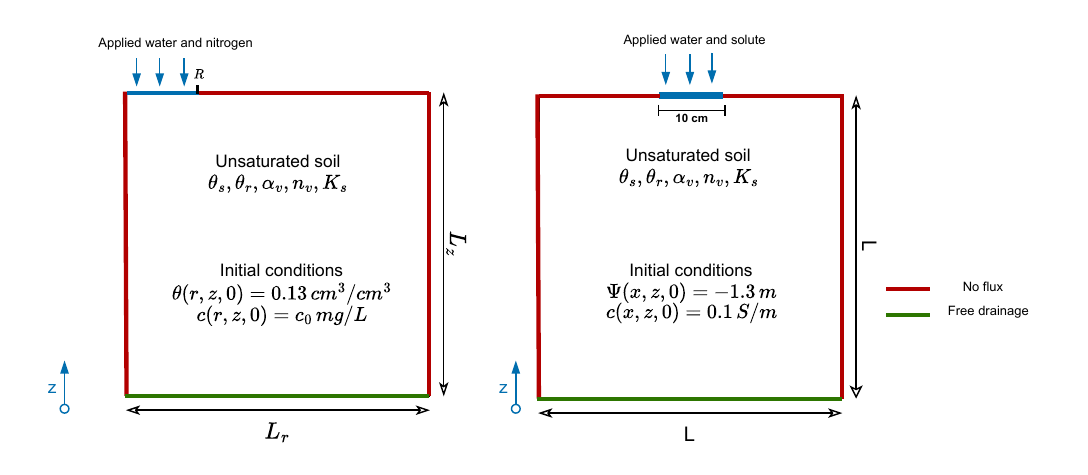}
\caption{Schematic illustration of the domain for the nitrate transport (left) and soil salt transport (right) tests.}
\label{Nitrate_PINNs_geometriy}
\end{center}
\end{figure}
\begin{figure}
\centering
\includegraphics[width=\linewidth, height=9cm]{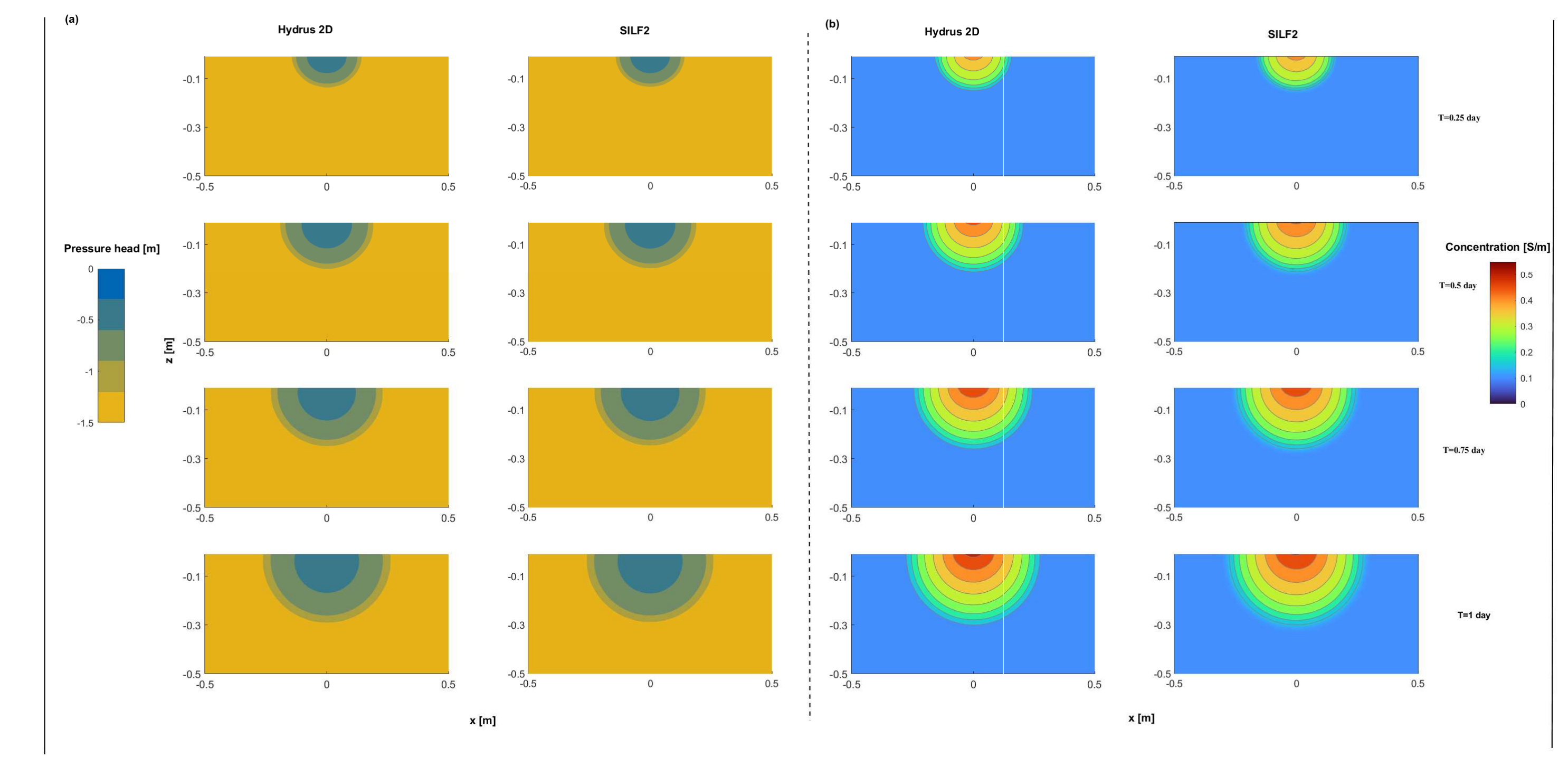}
\caption{Time evolution of the pressure head (a) and the salt soil transport (b) using SILF2 method compared to the Hydrus simulations in \cite{haruzi_moreno_2023}.}
\label{Fig:PINNS_Solute}
\end{figure}
\begin{figure}
\centering
\includegraphics[width=0.7\linewidth]{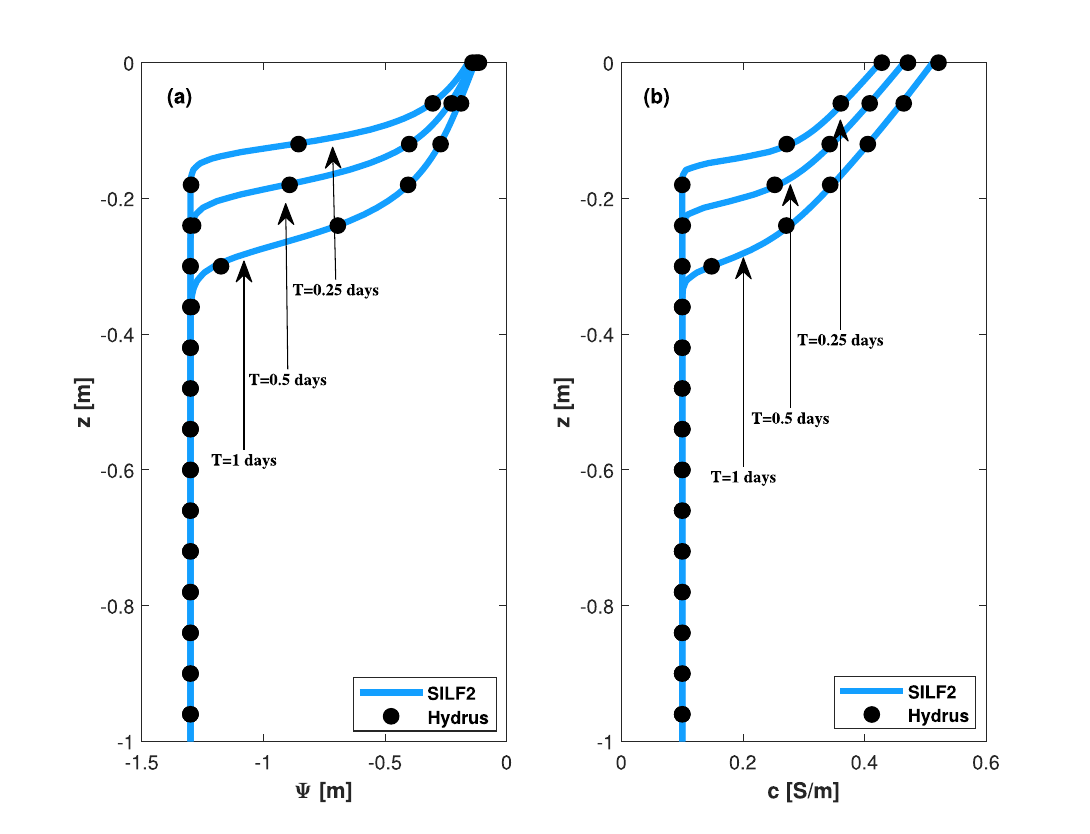}
\caption{ Pressure head (a) and soil salt transport (b) vertical cross-section ($x=0$). SILF2 scheme compared to the Hydrus simulations in \cite{haruzi_moreno_2023}.}
\label{Cross_PINNS}
\end{figure}

In this test, we will examine the proposed numerical approach for modeling water flow and soil salt transport through unsaturated porous media. Haruzi and Moreno \cite{haruzi2023modeling} investigate the infiltration and redistribution of water and soil salt transport in a loamy soil using a physics-informed neural network approach. As the authors described in \cite{haruzi2023modeling}, the soil salt concentration is determined by calculating the pore-water electrical conductivity (EC). The data used in this study is provided by the authors and it is available in \cite{haruzi_moreno_2023}. We will use the provided Hydrus simulations to compare them with the numerical results obtained using the proposed SILF2 scheme. A square computational domain with a length \(L = 2\,m\) is used. At the center of the top boundary, a section with a width of $10\, cm$ is considered to infiltrate the water using a rate of $0.1 \, m/day$. The initial pressure head condition is equal to \(-1.3\,m\). Free drainage is enforced at the bottom of the domain. The remaining boundaries are considered to have no water flux. Regarding the solute transport, the initial concentration condition is equal to $0.1\, S/m$. A third type of boundary condition (Cauchy boundary conditions) is applied in the section at the top boundary, involving a constant inflow of solute where the inflow of water is considered with constant concentration $c= 1\, S/m$. The VG model is used for the capillary pressure. Figure \ref{Nitrate_PINNs_geometriy} (right) summarizes the prescribed geometry and boundary conditions, while Table \ref{tab:pinns} provides parameters and material proprieties related to the test problem. The domain is discretized using an unstructured triangular mesh with a mesh size of 0.061 where a total of 19270 triangles, 9869 vertices, and 50 nodes on each side of the square domain are used. We consider a time step $\Delta t=10^{-3}\, day$ and the simulation time is $T=1\,day$. Figure \ref{Fig:PINNS_Solute}\textbf{a} and \ref{Fig:PINNS_Solute}\textbf{b} illustrate the evolution of the pressure head and solute transport with time, respectively. Good agreements are obtained between the results of the SILF2 method and those obtained in \cite{haruzi2023modeling}. For a more details comparison, we make a vertical cross-section at the center of the domain ($x=0$) for the pressure head and solute concentration solutions. The resulting 1D solutions are depicted in Figure \ref{Cross_PINNS} demonstrating good agreement between our results and those in \cite{haruzi2023modeling}.
\setlength{\arrayrulewidth}{1pt}
\noindent
\begin{table}[ht]
\centering

\begin{tabularx}{\textwidth}{>{\hsize=2\hsize} XX >{\hsize=0.4\hsize} XX}
\hline
\textbf{Material propriety} & \textbf{Symbol} & \textbf{value}  \\
\hline
Saturated hydraulic conductivity  & $K_s$ [$m /day$] & $0.25$ \\
Saturated water content & $\theta_s$ [$m^3/m^3$] \;\, & $0.43$  \\
Residual water content & $\theta_r$ [$m^3/m^3$] & $0.078$  \\
VG fitting coefficient   & \(\alpha_v \)  [$1/m$] & $3.6$ \\
VG parameter   & \(n_v\)  [$-$] & $1.56$ \\
Molecular diffusion coefficient & $\lambda_m$ [$m^2/day$] & $0.0$ \\
Longitudinal dispersivity & $\lambda_L$ [$m$] & $0.5$ \\
Transversal dispersivity & $\lambda_T$ [$m$] & $0.1$ \\

\hline
\end{tabularx}
\caption{Material proprieties used for the soil salt transport test \cite{haruzi2023modeling}. }
\label{tab:pinns}
\end{table}

%% file: Docts/Nitrate.tex
Further numerical investigation will be carried out to test the performance of the proposed methodology for modeling water flow and nitrate transport under surface drip fertigation. Li et al. 2005 \cite{li2005modeling} incorporate the presented test problem to compare numerical results with the laboratory data obtained from studies in \cite{li2003water}. The experiment was conducted in a cylinder where the infiltration and nitrate transport processes are considered axisymmetrical, emphasizing the radius ($r$) and depth ($z$) as variables in the analysis. Thanks to the axisymmetrical nature of the problem, a rectangular computational domain of length $L_r = 41 \, cm$ and width $L_z = 40\, cm$ is considered and a water entry zone of radius $R$ is located at the left of the top boundary as shown in Figure \ref{Nitrate_PINNs_geometriy} (left). In \cite{li2005modeling}, the authors investigate the impact of different fertigation strategies on water movement and nitrogen dynamics in the soil by applying ammonium nitrate fertilizer ($NH_4NO_3$) from the saturated entry zone. The VG model is used for the capillary pressure. The initial water content and nitrate concentration were distributed uniformly in the domain, that is, \(\theta(r, z, 0) = 0.13\, cm^3/ cm^{3}\) and \( c(r, z, 0)=247.5\, mg/L\). For the VG model, the capillary pressure is 0 in the saturated entry zone, and since the pounded water depth at the saturated entry zone is about $0.15 \, cm$, the hydrostatic pressure will be added, i.e. $\Psi = 0.15 \, cm$. The boundary conditions for this test are expressed as follows:
\begin{itemize}
 \item Dirichlet boundary conditions for the pressure head and nitrate concentration are imposed at the surface saturated entry zone:
\begin{equation}
\left\{
    \begin{array}{ll}
\begin{aligned}
\Psi(r,L_z,t)& =0.15 \; cm & 0\leq r \leq R, \\
c(r,L_z,t)& =C_a \; mg/L & 0\leq r \leq R.     
\end{aligned}
    \end{array}
    \right.
\end{equation}
\item No flux boundary conditions  for water flow are applied to the top boundary of the domain outside of the surface saturated entry zone:
$$
-K\Big[\dfrac{\partial\Psi}{\partial r} + 1 \Big]=0 \quad R \leq r \leq L_r, z=L_z.
$$
\item No-flux boundary conditions for water flow and nitrate concentration are enforced in the lateral sides of the domain:
\begin{equation}
\left\{
    \begin{array}{ll}
\begin{aligned}
-K\dfrac{\partial \Psi}{\partial r}=0 \quad  & r=0, \; r=L_r, \; 0\leq z \leq L_z, \\[5pt]
\theta D_{rr}\dfrac{\partial c}{\partial r}=0 \quad  & r=0, \; r=L_r, \; 0\leq z \leq L_z.     
\end{aligned}
    \end{array}
    \right.
\end{equation}
\item Free drainage boundary conditions are used in the bottom of the domain:
\begin{equation}
\left\{
    \begin{array}{ll}
\begin{aligned}
\dfrac{\partial \Psi}{\partial z}=0 \quad  & z=0, \; 0\leq r \leq L_r, \\[5pt]
\theta D_{zz}\dfrac{\partial c}{\partial z}=0 \quad  & z=0, \; 0\leq r \leq L_r.  \qquad \qquad \quad    
\end{aligned}
    \end{array}
    \right.
\end{equation}
\end{itemize}

\begin{figure}
\centering
\includegraphics[width=0.65\linewidth]{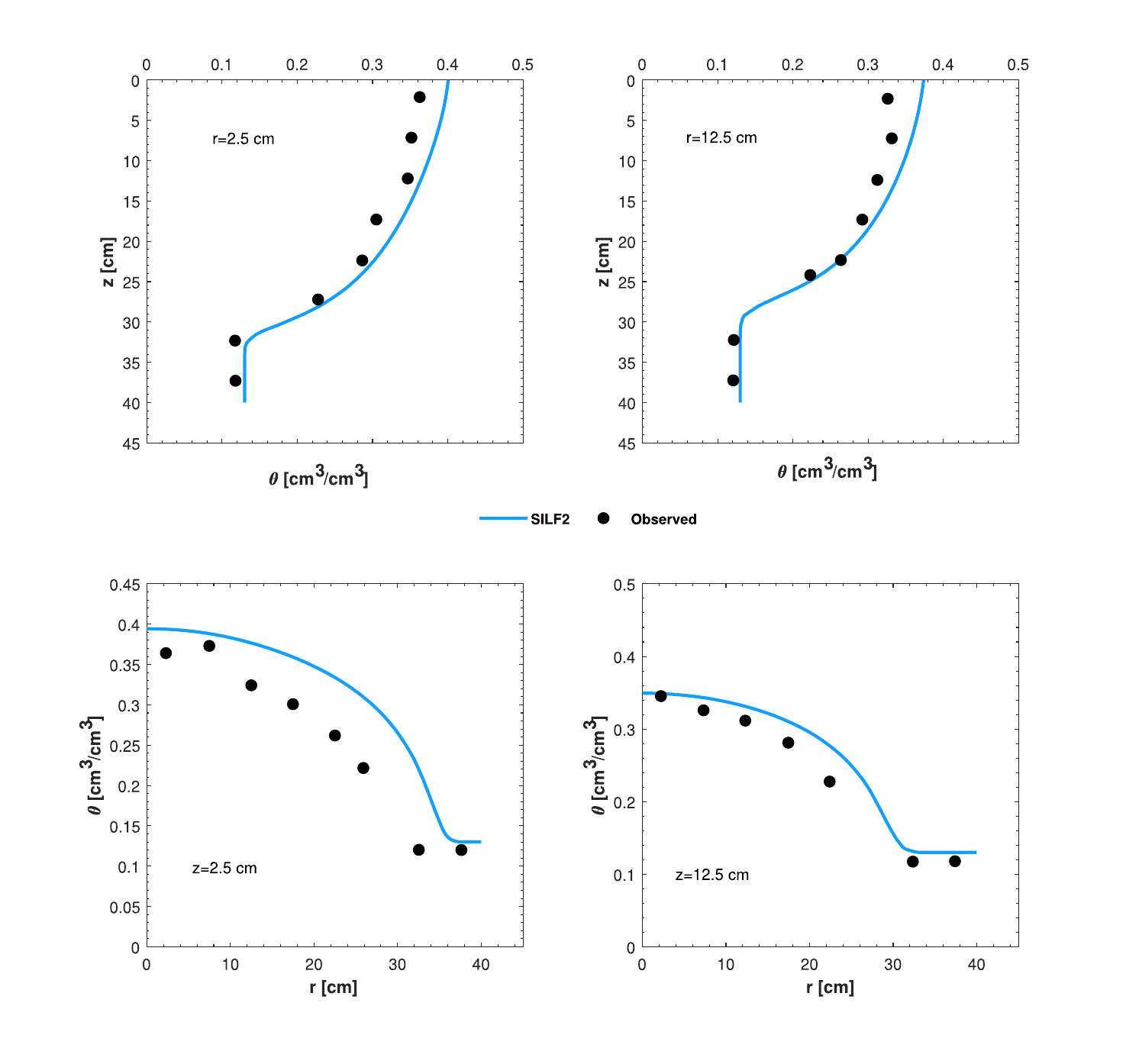} 
\caption{Comparison of the observed and simulated water content for the nitrate transport test.}
\label{Li_WC}
\end{figure}

\begin{figure}
\centering
\includegraphics[width=0.8\linewidth,height=10cm]{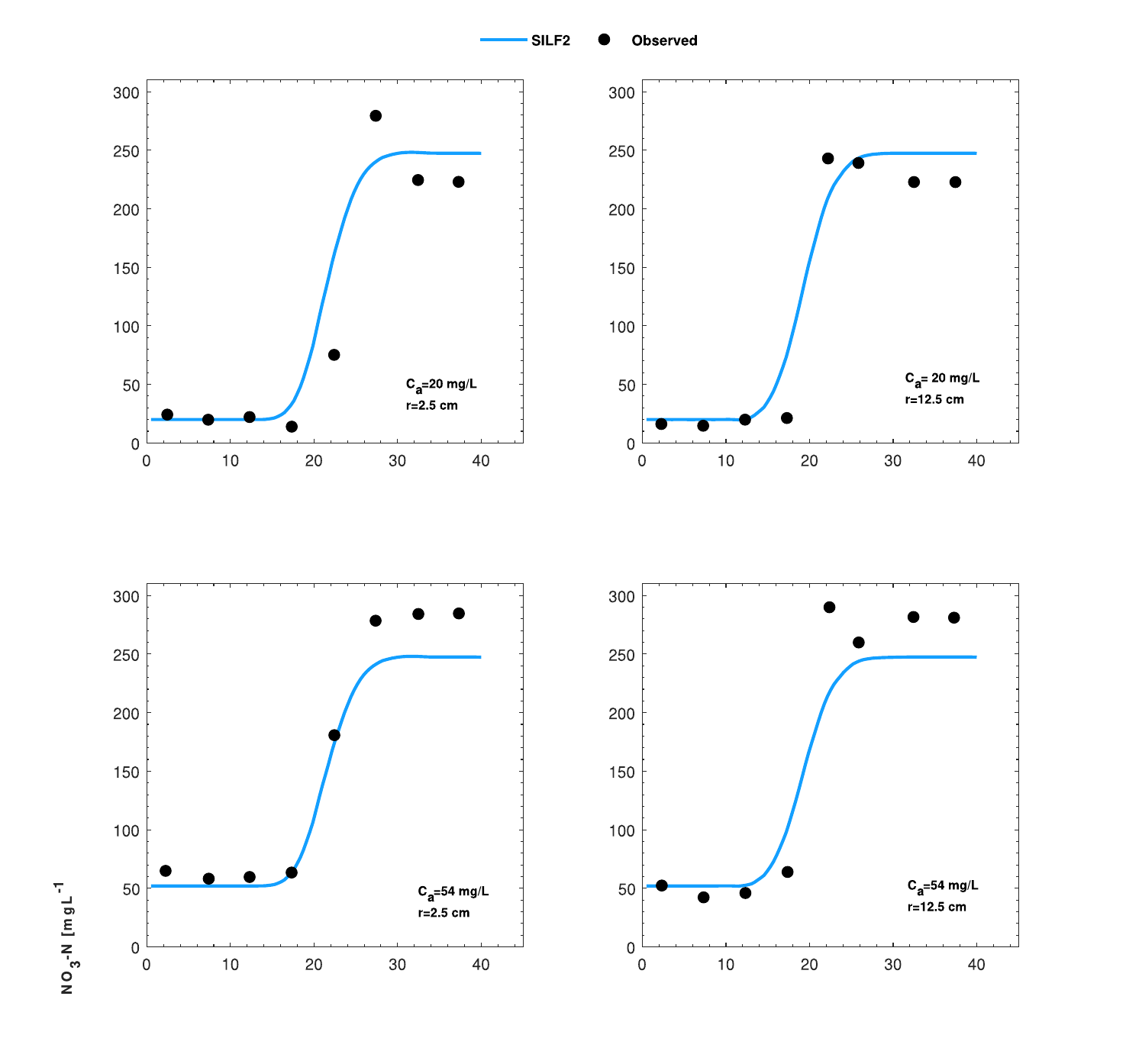} 
\vspace{0cm}
\includegraphics[width=0.8\linewidth, height=10cm]{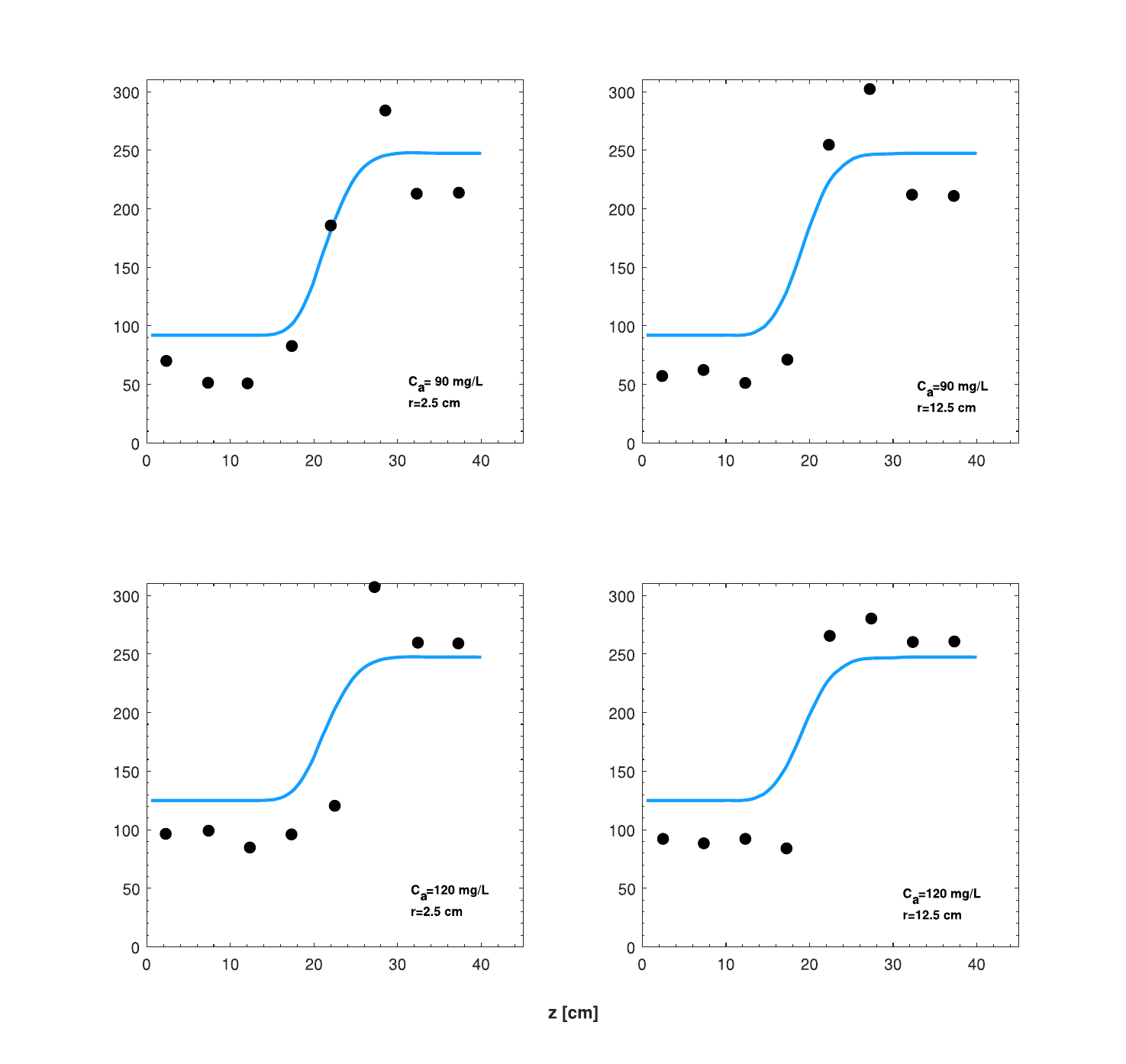}
\caption{Comparison between the simulated and observed results of the nitrate concentration using different inputs of nitrate concentration.}
\label{Li_Nitrate}
\end{figure}

\setlength{\arrayrulewidth}{1pt}
\noindent
\begin{table}[ht]
\centering

\begin{tabularx}{\textwidth}{>{\hsize=2\hsize} XX >{\hsize=0.4\hsize} XX}
\hline
\textbf{Material propriety} & \textbf{Symbol} & \textbf{value}  \\
\hline
Radius of the entry zone  & $R$ [$cm$] & $5.6$ \\
Saturated hydraulic conductivity  & $K_s$ [$cm/hour$] & $1.96$ \\
Saturated water content & $\theta_s$ [$cm^3/cm^3$] \;\, & $0.41$  \\
Residual water content & $\theta_r$ [$cm^3/cm^3$] & $0.047$  \\
VG fitting coefficient   & \(\alpha_v \)  [$1/cm$] & $0.015$ \\
VG parameter   & \(n_v\)  [$-$] & $1.48$ \\
Molecular diffusion coefficient & $\lambda_m$ [$cm^2/min$] & $0.0015$ \\
Longitudinal dispersivity & $\lambda_L$ [$cm$] & $0.32$ \\
Transversal dispersivity & $\lambda_L$ [$cm$] & $0.0032$ \\

\hline
\end{tabularx}
\caption{Parameters and material proprieties used for nitrate transport under surface drip fertigation \cite{li2005modeling}.}
\label{tab:nitrate}
\end{table}\\
Parameters and material properties related to the test problem are listed in Table \ref{tab:nitrate}. In \cite{li2005modeling}, \(NH_4NO_3\) fertilizer is applied, where mineralization and some reactions are neglected. The adsorption process between the solid phase of the soil and nitrate is also neglected. Regarding the numerical method, the domain geometry is discretized by an unstructured mesh with a total of 3153 non-uniform triangles and 1651 vertices and a time step  \(\Delta t= 10^{-4}\) is used. The water movement and nitrate transport in the domain are simulated at the final time $T=8$ hours. Figure \ref{Li_WC} depicts the comparison of the simulated and experimental water content using cross-section profiles. Vertical cross-sections are fixed at $r=2.5\, cm$ and $r=12.5\, cm$, while  horizontal cross-sections are fixed at $z=2.5\, cm$ and $z=12.5\, cm$. Furthermore,  Figure \ref{Li_Nitrate} illustrates vertical cross-section profiles, fixed at $r=2.5\, cm$ and $r=12.5\, cm$, of the nitrate concentration for various input concentrations $C_a$ of  \(NO_3^-\)  as a function of depth. The numerical results obtained using the SILF2 method are in good agreement with the observed data.

\begin{figure}
\centering
\includegraphics[width=\linewidth]{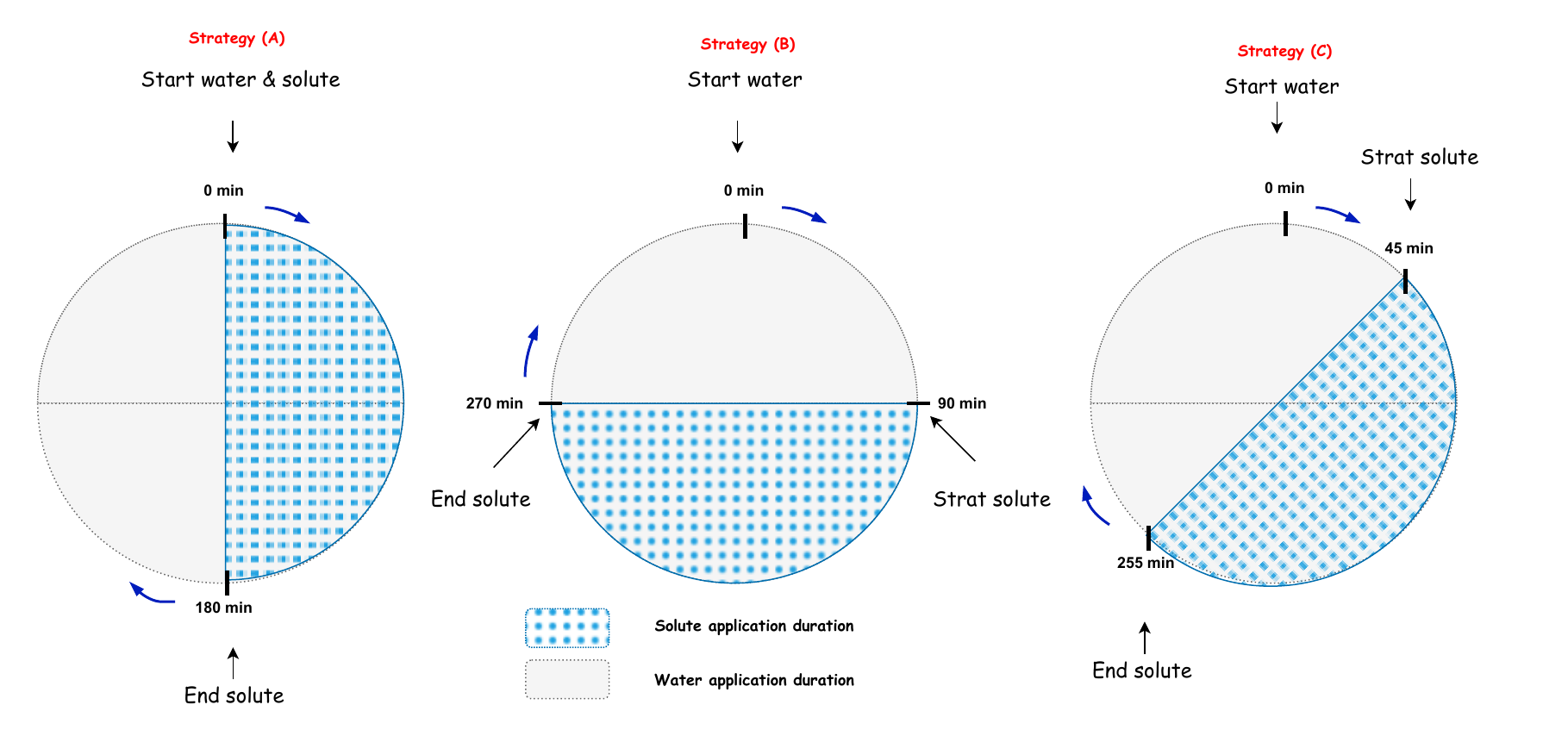} 
\caption{Visual representations of fertigation strategies. The total time irrigation is 360 min.}
\label{Fig:Strategies}
\end{figure}

\begin{figure}
\centering
\includegraphics[width=\linewidth]{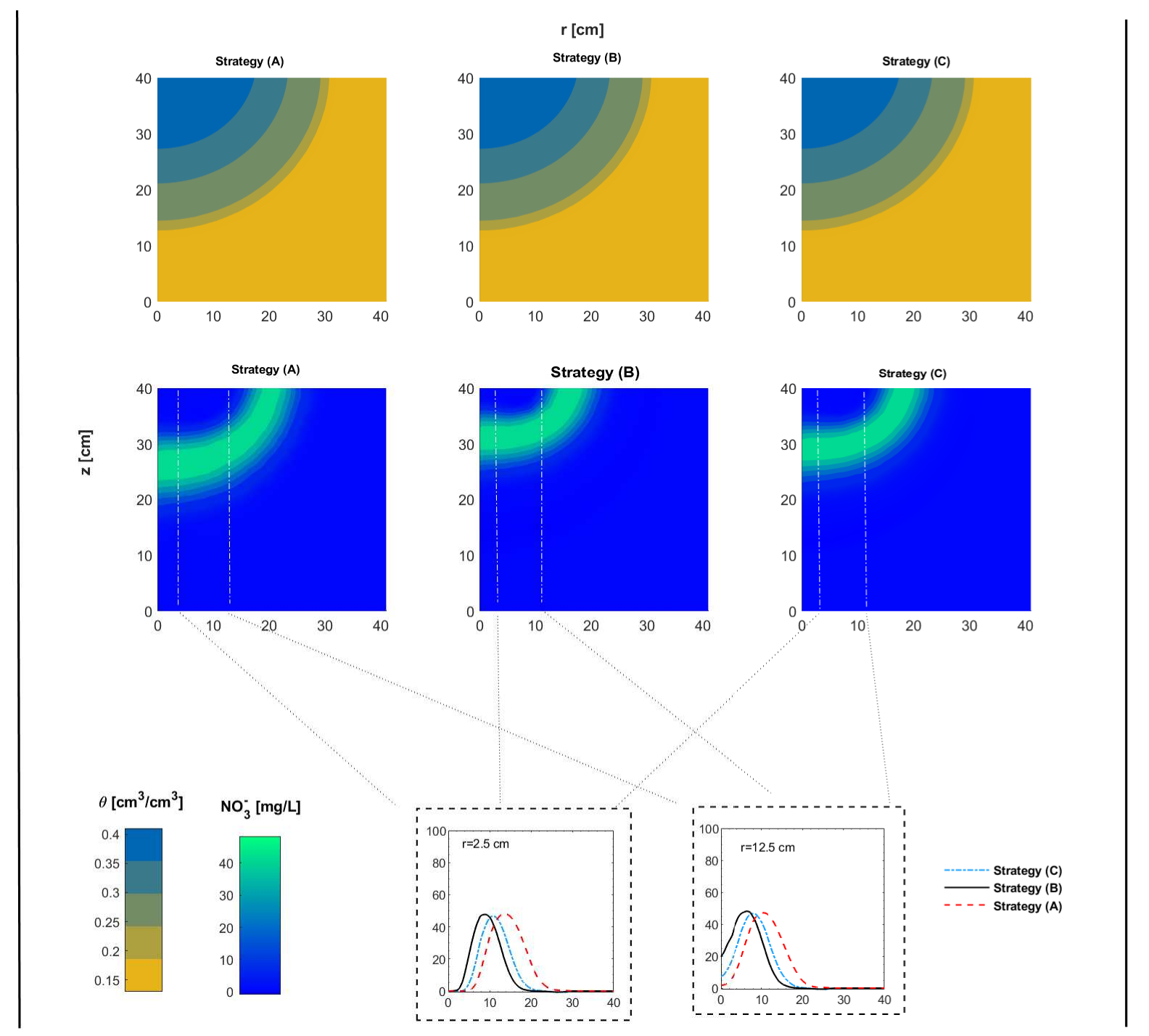} 
\caption{Water content (top) and nitrate (down) distributions after 6 hours of fertigation for the three strategies using the SILF2 method.}
\label{Li_Stratigies}
\end{figure}

In \cite{li2005modeling}, the authors suggested three different strategies for drip fertigation of nitrate in the soil to evaluate the water and nitrate use efficiency:
\begin{itemize}
    \item Strategy (A): Apply water and \(NH_4NO_3\) at the beginning for $1/2$ of the total irrigation duration, followed by continued application of water for the remaining time of irrigation.
    \item Strategy (B): Apply water, firstly, for $1/4$ of the total irrigation duration, then apply \(NH_4NO_3\) for $1/2$ of the total irrigation duration, followed by continued application of water for the remaining $1/4$ of the total time irrigation.
    \item Strategy (C): Apply water, firstly, for $1/8$ of the total irrigation time, then apply \(NH_4NO_3\) for $1/2$ of the total irrigation duration, followed by continued application of water for the remaining $3/4$ of the total time irrigation. Figure \ref{Fig:Strategies} shows a visual depiction of these fertigation strategies.
\end{itemize}
Throughout the studies, the authors demonstrate that Strategy (B) is recommended as the optimal fertigation strategy. We will evaluate the proposed SILF2 method in addressing the simulation of these three strategies. The discretization of the domain, the time step, and the initial water content are preserved as in the previous test, while the initial concentration of the nitrate is assigned the value  $0\, mg/L$. The spatial distribution of water content and nitrate concentration for the three strategies at the time $T=6$ hours are illustrated in Figure \ref{Li_Stratigies}. Additionally, Figure \ref{Li_Stratigies} depicts vertical cross-section profiles, with $r=2.5\, cm$ and $r=12.5\, cm$, of nitrate concentration. We obtained good agreement between the numerical results of the SILF2 method and those presented in \cite{li2005modeling}. 

%% file: Docts/TError.tex
To demonstrate the second-order convergence of the proposed numerical schemes, we consider the following equation:
\begin{equation}
    \dfrac{\partial \varphi}{\partial t} = F\left(\varphi, t \right),
    \label{A1}
\end{equation}
with the following temporal discretization:
\begin{equation}
  \dfrac{[\![ \varphi(t_{n+\delta}) ]\!]}{\Delta t} = \Big[A(\varphi(t_{n+\delta})\Big],
  \label{A2}
\end{equation}
where
\begin{equation}
 \begin{aligned}
\left[\!\left[ \varphi(t_{n+\delta}) \right]\!\right] &= \left(\frac{1}{2} + \delta\right)\varphi\left(t_{n+1}\right)-2\delta\varphi(t_{n})+\left(\delta - \frac{1}{2}\right)\varphi(t_{n-1}),\\ \Big[A(\varphi(t_{n+\delta})\Big]&=\left(\delta+\mu\right)F\big(\varphi(t_{n+1})\big)+ \left(1-\delta - 2\mu\right) F\big(\varphi(t_n)\big) + \mu F\big(\varphi(t_{n-1})\big).
\end{aligned}
\label{A.3}
\end{equation}
By employing the Taylor series expansion and with the assumption of sufficient smoothness, the terms within equation \eqref{A.3} can be centered about the current time level \( t_{n+\delta}\) as follows:


\begin{equation}
    \begin{aligned}
    &\dfrac{[\![ \varphi(t_{n+\delta}) ]\!]}{\Delta t}=\dfrac{\partial\varphi}{\partial t}(t_{n+\delta}) + \frac{\Delta t^2}{6} \left[ (1 - \delta)^3 \Big( \delta + \frac{1}{2} \Big)  + 2\delta^4  - (1 + \delta)^3 \Big( \delta - \frac{1}{2} \Big) \right] \frac{\partial^3 \varphi}{\partial t^3}(t_{n+\delta}) + O(\Delta t^3),\\
    &\Big[A(\varphi(t_{n+\delta})\Big] = F(\varphi(t_{n+\delta}))+ (\delta +\mu) \bigg[\Big(\varphi(t_{n+1})-\varphi(t_{n+\delta})\Big)F'(\varphi(t_{n+\delta})) + \dfrac{\left(\varphi(t_{n+1})- \varphi(t_{n+\delta})\right)^2}{2}F''(\varphi(t_{n+\delta}))\bigg]\\  
    &\qquad \qquad \; \quad + (1-\delta -2\mu)\bigg[\Big(\varphi(t_{n})-\varphi(t_{n+\delta})\Big)F'(\varphi(t_{n+\delta})) + \dfrac{(\varphi(t_{n})- \varphi(t_{n+\delta}))^2}{2}F''(\varphi(t_{n+\delta}))\bigg]\\
     &\qquad \qquad \; \quad +\mu\bigg[\Big(\varphi(t_{n-1})-\varphi(t_{n+\delta})\Big) F'(\varphi(t_{n+\delta})) + \dfrac{(\varphi(t_{n-1})- \varphi(t_{n+\delta}))^2}{2}F''(\varphi(t_{n+\delta}))\bigg]\\
    &\qquad \qquad \; \quad  +O \left(\Delta_{n+1-\delta}^3\right)+O\left(\Delta_{n+\delta}^3\right)+O\left(\Delta_{n-1-\delta}^3\right),
    \end{aligned}
    \label{A4}
\end{equation}

where $\Delta_{n+1-\delta} = \varphi(t_{n+1})-\varphi(t_{n+\delta})$,  
    $\Delta_{n+\delta}= \varphi(t_{n})-\varphi(t_{n+\delta})$ and 
    $\Delta_{n-1-\delta}= \varphi(t_{n-1})-\varphi(t_{n+\delta})$.\\
    
Using the Taylor series, the term \(\varphi(t_{n+1})\) is approximated around the time level $t_{n+\delta}$ as follows:
\begin{equation}
    \begin{aligned}
    &\varphi(t_{n+1})= \varphi(t_{n+\delta})+\Delta t(1-\delta)\dfrac{\partial \varphi}{\partial t}\left(t_{n+\delta}\right)+\dfrac{\Delta t^2 (1-\delta)^2}{2}\dfrac{\partial^2 \varphi}{\partial t^2}\left(t_{n+\delta}\right) +O(\Delta t^3).
\end{aligned}    
\end{equation}
By utilizing the same calculations to expand \(\varphi(t_{n})\) and \(\varphi(t_{n-1})\) around the time level $t_{n+\delta}$, we get
\begin{equation}
\begin{aligned}
\Big[A(\varphi(t_{n+\delta})\Big] = F(\varphi(t_{n+\delta})) &+ \Delta t^2 \bigg[ \dfrac{(\delta + \mu)(1-\delta)^2}{2} + \dfrac{ (1-\delta - 2\mu)\delta^2 }{2} + \dfrac{\mu(1+\delta)^2}{2}\bigg]\\
&\times\Big( F'(\varphi(t_{n+\delta})) + F''(\varphi(t_{n+\delta})) \Big) + O(\Delta t^3).
\end{aligned}
\end{equation}
We conclude that the scheme \eqref{A2} exhibits a second-order truncation error. Consequently, BDF2, SBDF2 and CN2 schemes are of second order accuracy.

 Regarding the SILF2 method presented in \eqref{SILf}, we will use the following notations:
 \begin{equation}
 \begin{aligned}
 [\![ \Psi_h(t_n) ]\!] &= \dfrac{ \Psi_h(t_{n+1})-\Psi_h(t_{n-1})}{2\Delta t},\\
[\![ R_h(t_n) ]\!]_\nu &= \nu \Psi_h(t_{n+1}) +(1-2\nu)\Psi_h(t_{n})+\nu \Psi_h(t_{n-1}).   
 \end{aligned}
 \label{A.7}
\end{equation}

By employing the Taylor series to expand the terms in \eqref{A.7}, we obtain:
\begin{equation}
\begin{aligned}
     &[\![ \Psi_h(t_n) ]\!]= \dfrac{\partial\Psi_h}{\partial t}(t_{n}) + \left(\dfrac{\Delta t^2}{6}\right) \frac{\partial^3 \Psi_h}{\partial t^3}(t_{n}) + O(\Delta t^3),\\
    &\qquad \qquad \;  = \dfrac{\partial\Psi}{\partial t}(t_{n}) + O(\Delta t^2),\\
    & [\![ R_h(t_n) ]\!]_\nu = \Psi_h(t_n) + \Delta t^2 \frac{\partial^2 \Psi_h}{\partial t^2}(t_{n}) + O(\Delta t^3).
\end{aligned}
\end{equation}

We conclude that the SILF2 method exhibits a second-order truncation error.